\documentclass{amsart}

\usepackage{amssymb}

\def\bc{\begin{center}}
\def\ec{\end{center}}
\def\be{\begin{equation}}
\def\ee{\end{equation}}
\def\ben{\begin{enumerate}}
\def\een{\end{enumerate}}
\def\bfg{\begin{figure}}
\def\efg{\end{figure}}
\def\bq{\begin{quote}}
\def\eq{\end{quote}}
\def\bd{\begin{description}}
\def\ed{\end{description}}
\def\h{\hbar}
\def\p{\partial}
\def\w{\wedge}

\def\dim{\operatorname{dim}}

\def\tr{\operatorname{tr}}

\newcommand{\CC}{{\mathbb C}}

\newcommand{\ZZ}{{\mathbb Z}}

\newcommand{\lan}{\langle}
\newcommand{\ran}{\rangle}
\newcommand{\Gb}{{\mathfrak B}}

\newcommand{\ga}{\alpha}

\renewcommand{\ge}{\varepsilon}

\newcommand{\gl}{\lambda}

\newcommand{\go}{\omega}

\newcommand{\gs}{\sigma}

\newcommand{\gO}{\Omega}

\newcommand{\J}{{\mathcal J}}

\newcommand{\M}{\overline{\mathcal M}}

\renewcommand{\CC}{\mathbb C}
\renewcommand{\t}{\mathbf t}
\newcommand{\q}{\mathbf q}
\newcommand{\f}{\mathbf f}
\newcommand{\g}{\mathbf g}
\newcommand{\m}{\mathbf m}
\newcommand{\x}{\mathbf x}
\newcommand{\y}{\mathbf y}

\newcommand{\D}{\mathcal D}

\newcommand{\A}{\mathcal A}
\newcommand{\W}{\mathcal W}
\renewcommand{\a}{\alpha}
\renewcommand{\b}{\beta}
\renewcommand{\c}{\gamma}
\renewcommand{\d}{\delta}
\renewcommand{\H}{{\mathcal H}}
\newcommand{\F}{{\mathcal F}}

\newcommand{\T}{{\mathcal T}}
\newcommand{\C}{{\mathcal C}}
\renewcommand{\L}{{\mathcal L}}

\newcommand{\1}{{\bf 1}}

\newcommand{\0}{{\mathbf 0}}
\newcommand{\Res}{\operatorname{Res}}

\newcommand{\Beta}{\Gb}
\newcommand{\diag}{\operatorname{diag}}

\begin{document}

\title[Simple singularities and integrable hierarchies ]
{ Simple singularities and integrable hierarchies}
%%%% \dedicatory{} 
\author{Alexander B. Givental \and Todor E. Milanov} 
\address{UC Berkeley} 
\thanks{Research is partially supported by NSF Grant DMS-0072658} 

%%%% \date{ }

\begin{abstract}

The paper \cite{GiQ} gives a construction of the 
{\em total descendent potential} corresponding to a semisimple Frobenius 
manifold. In \cite{GiI}, it is proved   
that the total descendent potential corresponding to K. Saito's
Frobenius structure on the parameter space of the miniversal deformation of the
$A_{n-1}$-singularity satisfies the modulo-$n$ reduction of the KP-hierarchy.
In this paper, we identify the hierarchy satisfied by the total
descendent potential of a simple singularity of the $A,D,E$-type. 
Our description of the hierarchy is parallel to the vertex operator 
construction of Kac -- Wakimoto \cite{KW} except that we give both some 
general integral formulas and explicit numerical values for certain 
coefficients which in the Kac -- Wakimoto theory are studied on 
a case-by-case basis and remain, generally speaking, unknown.

\end{abstract}

\maketitle

%\section{Introduction}

\subsection*{\bf 1. The ADE-hierarchies.}
The KdV-hierarchy of integrable systems can be placed under the name $A_1$
into the list of more general integrable hierarchies corresponding to 
the $ADE$ Dynkin diagrams. These hierarchies are usually constructed \cite{K}
using representation theory of the corresponding loop groups. 
V. Kac and M. Wakimoto \cite{KW} describe the hierarchies even more
explicitly in the form of the so called {\em Hirota quadratic equations} 
expressed in terms of suitable {\em vertex operators}.

One of the goals of the present paper is to show how the vertex operator 
description of the Hirota quadratic equations (certainly the same ones, 
even though we don't quite prove this) emerges from the theory of vanishing
cycles associated with the $ADE$ singularities. 

Let $f$ be a weighted-homogeneous polynomial in $\CC^3$ with a simple 
critical point at the origin. According to V. Arnold \cite{Ar} simple
singularities of holomorphic functions are classified by 
simply-laced Dynkin diagrams:
\[ A_{N},\ N\geq 1: \ \ f=\frac{x_1^{N+1}}{N+1}+\frac{x_2^2}{2}+\frac{x_3^2}{2}
, \ \ 
D_N,\ N\geq 4: \ \ f=x_1^2x_2-x_2^{N-1}+x_3^2, \]
\[ E_6: \ \ f=x_1^3+x_2^4+x_3^2,\ \ E_7:\ \ f=x_1^3+x_1x_2^3+x_3^2,\ \ 
E_8: \ \ f=x_1^3+x_2^5+x_3^2. \]
Let $H=\CC [x_1,x_2,x_3]/(f_{x_1},f_{x_2},f_{x_3})$ denote the local algebra 
of the critical point. We equip $H$ with a non-degenerate symmetric 
bilinear form $(\cdot,\cdot )$ by picking a weighted - homogeneous 
holomorphic volume $\go = dx_1\w dx_2\w dx_3$ 
and using the residue pairing:
\[ ( \varphi, \psi ) :=\operatorname{Res}_{0} \frac{\varphi(x)\ \psi(x)\ \go}
{f_{x_1}f_{x_2}f_{x_3}}.\]
Let $\H = H ((z^{-1}))$ be the space of Laurent series 
$\f(z)=\sum_{k\in \ZZ} \f_k z^k$ in one indeterminate 
$z^{-1}$ ({\em i.e.} finite in the direction of positive $k$) 
with vector coefficients $\f_k\in H$. 
We endow $\H$ with the symplectic form
\[ \gO (\f,\g):= \frac{1}{2\pi i} \oint ( \f(-z),\g(z) )\ dz .\] 
The polarization $\H = \H_{+}\oplus \H_{-}$ where 
$\H_{+}=H[z]$ and $\H_{-}=z^{-1}H[[z^{-1}]]$ is Lagrangian and identifies 
$\H$ with the cotangent bundle space $T^*\H_{+}$. 
The Hirota quadratic equations are imposed on {\em asymptotical functions} of 
$\q = \q_0+\q_1z+\q_2z^2+...\ \in \H_{+}$. By an asymptotical function
we mean an expression of the form
\[ \Phi = \exp \sum_{g=0}^{\infty} \h^{g-1} \F^{(g)}(\q ) \]
where usually $\F^{(g)}$ will be formal functions on $\H_{+}$.
By definition, vertex operators are elements of the Heisenberg group
acting on such functions. Given a sum $\f = \sum \f_k z^k $ (possibly infinite
in both directions) one defines the corresponding vertex operator of the form
\[ e^{ \Omega (\f_{-}, \q ) /\sqrt{\h} } e^{\sqrt{\h} \p_{\f_{+}}} =
\exp \{ \sum_{k\geq 0} (-1)^{k+1} \sum_{a}f^{a}_{-1-k}q_k^{a}/\sqrt{\h} \} 
\ \exp \{ \sqrt{\h} \sum_{k\geq 0}\sum_{a} f_k^{a} \p/\p q_k^{a} \}.\]
Here $f^{a}_k, q^{a}_k$ are components of the vectors $\f_k, \q_k$ in 
an {\em orthonormal} basis. 
We will make use of the vertex operators $\Gamma^{\phi}(\gl)$
corresponding to $2$-dimensional homology classes 
$\phi \in H_2(f^{-1}(1)) \simeq \ZZ^N$ and defined 
as follows. Take 
\[ \f := \sum_{k\in \ZZ} I^{(k)}_{\phi} (\gl) (-z)^k\ \ \text{where}\ \ 
d I^{(k)}_{\phi}/d\gl = I^{(k+1)}_{\phi}, \]
and $I^{(-1)}_{\phi}(\gl) \in H$ is the following period vector:
\[ ( I^{(-1)}_{\phi}(\gl) , [\psi_{a} ] ) := \frac{1}{2\pi}
\int_{\phi \subset f^{-1}(\gl)} \psi_{a}(x) \frac{\go}{df} .\]
The cycle $\phi$ is transported from the level surface $f^{-1}(1)$ to 
$f^{-1}(\gl)$, and $\psi_{a}$ are {\em weighted-homogeneous}
functions representing a basis in the local algebra $H$.
\footnote{As it follows, for instance, from \cite{GiA}, 
the integral on the R.H.S. depends only on the class $[\psi_a] \in H$.}
The functions 
$(I_{\phi}^{(k)}, [\psi_{a} ])$ are proportional to the fractional powers 
$\gl^{m_{a}/h-k-1}$ where $h$ is the Coxeter number and 
$m_{a}=1+h \deg \psi_{a}$ are the exponents of the appropriate 
reflection group $A_N,D_N$ or $E_N$. 

The lattice $H_2(f^{-1}(1))$ carries the action of the monodromy group
(defined via morsification of the function $f$) which is the reflection 
group with respect to the intersection form of cycles. The form is negative
definite, and we will denote $\lan \cdot, \cdot\ran$ the positive
definite form opposite to it. Let $A$ denote the set of 
{\em vanishing cycles},
{\em i.e.} the set of classes $\a \in H_2(V_1)$ with $\lan \a,\a \ran =2$ 
such that the reflections $\phi \mapsto \phi - \lan \a, \phi\ran \a$ 
belong to the monodromy group. The {\em Hirota quadratic equation} 
of the ADE-type takes on the form
\begin{align} \label{Hirota}
\Res_{\gl=\infty}\ \frac{d\gl}{\gl}\ \left[ \sum_{\a\in A} 
a_{\a}\ \Gamma^{\a}(\gl) \otimes \Gamma^{-\a}(\gl) \right] \ 
(\Phi \otimes \Phi)   =  
\frac{N (h+1)}{12h} \ (\Phi \otimes \Phi) + \\ 
\label{L-nod} 
 \sum_{k\geq 0} \sum_a (\frac{m_a}{h}+k) (q^a_k\otimes 1 -1\otimes q^a_k)
(\frac{\p}{\p q^a_k}\otimes 1-1\otimes \frac{\p}{\p q^a_k}) \     
(\Phi \otimes \Phi). \end{align}
The tensor product sign means that the functions depend on two copies 
$\q'$ and $\q''$ of the variable $\q$, 
and the objects on the left of $\otimes$
refer to $\q=\q'$ while those on the right --- to $\q=\q''$. 
The equation can be interpreted as follows. Put $\q'=\x+\y, \q''=\x-\y$.
and expand (\ref{Hirota},\ref{L-nod}) as a power series in $\y$.
Namely, rewrite the vertex operators: 
\[ \Gamma^{\a}(\gl)\otimes \Gamma^{-\a}(\gl)= \exp \{
\sum 2 (-1)^{k+1}f_{-1-k}^a  \h^{-1/2}y_k^a\}\ \exp \{
\sum f_k^a \h^{1/2} \p_{y_k^a}\} \]
where the coefficients $f_k^a$ (respectively $f_{-1-k}^a$) are proportional 
to negative (respectively positive) fractional powers of $\gl$. 
The residue sum (which should be understood here as the coefficient at 
$\gl^0$) can therefore be written as a power series 
$\sum\y^{\m} P_{\m}(\p_{\y})$ in $\y$ with coefficients $P_{\m}$ which 
are differential polynomials. Also, $\Phi \otimes \Phi = 
\Phi (\x+\y) \Phi (\x-\y)$ can be expanded into the Taylor power series in 
$\y$ with coefficients being quadratic expressions in partial derivatives of 
$\Phi (\x)$. Finally (\ref{L-nod}) assumes the form
$ 2 \sum_{a,k}\ (m_a/h\ +\ k) y_k^a \p_{y_k^a}$.
Equating coefficients in (\ref{Hirota},\ref{L-nod}) at the same monomials 
$\y^{\m}$ we obtain a hierarchy of quadratic relations between partial
derivatives of $\Phi(\x)$. 

In particular, the equation corresponding to $\y^{\0}$ shows that
\begin{equation} \label{c0} \sum_{\a\in A} a_{\a} = \frac{N (h+1)}{12 h} 
\end{equation}
is a necessary condition for consistency of the hierarchy ({\em i.e.}
for existence of a non-zero solution $\Phi$).

According to C. Hertling (see the last chapter in \cite{H}) for any 
weighted - homogeneous singularity the expressions $N (h+1)/(12 h)$ 
and $h^{-2}\sum_a m_a (h-m_a)/2$ coincide. 
Therefore the operator on the R.H.S. of the Hirota equation 
is twice the {\em Virasoro operator}\footnote{In a sense it corresponds 
to the vector field $\gl \p_{\gl}$ in the Lie algebra of vector fields on 
the line --- see Section $7$ for further information about this.}
 \[ \sum_{a,k}\ (\frac{m_a}{h}+k)\ y_k^a\ \p_{y_k^a}\ +\ 
\sum_a\ \frac{m_a (h-m_a)}{4h^2} .\]

The coefficients
$a_{\a}$ actually depend only on the orbit of the vanishing cycle $\a$
under the action of the {\em classical monodromy operator} 
defined by transporting the cycles in $f^{-1}({\gl})$ 
around $\gl =0$ and acting as one of the Coxeter elements 
in the reflection group. In fact the root system $A$ consists of $N$ such 
orbits with $h$ elements each. Summing the vertex operators within 
the same orbit acts as taking the average over all $h$ branches of the
function $\gl^{1/h}$. Thus the total sum does not contain fractional powers
of $\gl $ when expanded near $\gl =\infty$.   

The exact values of the coefficients $a_{\a}$ can be described as follows.
To a vector $\b \in H_2(f^{-1}(1),\CC)\simeq \CC^N$, associate the 
meromorphic $1$-form on $\CC^N$
\begin{equation} \label{W} 
\W_{\b}:= -\frac{1}{2}\sum_{\c \in A} \lan \b, \c \ran^2 
\frac{d \lan \c , x\ran}{\lan \c, x\ran } .\end{equation}
Let $w$ be an element of the reflection group and $\a$ and $\b=w\a$ be two
roots. Then
\begin{equation} \label{c} a_{\b} / a_{\a} = \exp \int_{\kappa}^{w^{-1}\kappa} \W_{\a} = \prod_{\c\in A} \lan \kappa,\c \ran^{\lan \a,\c\ran^2/2-
\lan \b,\c\ran^2/2}, \end{equation}  
where $\kappa \in \CC^N$ denotes an eigenvector of the classical 
monodromy operator $M$ with the eigenvalue $\exp (2\pi i/h )$. 
The R.H.S. does not depend on the 
path connecting $\kappa$ with $w^{-1}\kappa$ since $\W_{\a}$ is closed with 
logarithmic poles on some mirrors and with periods which are integer multiples 
of $2\pi i$. It does not depend on the normalization of $\kappa$ since 
$\W_{\a}$ is homogeneous of degree $0$. Also, the identity (see {\em i.g.} 
\cite{Bo}, Section $V.6.2$)  
\[ \sum_{\c\in A} \lan \c , x\ran^2 = 2h\ \lan x,x\ran \] 
implies that $i_{\sum x_a\p/\p x_a} \W_{\a} = -2h$ and shows that 
$\int_{\kappa}^{M^{-1}\kappa}\W_{\a} = h^{-1} \int_{\kappa}^{M^{-h}\kappa}
\W_{\a} = -4\pi i$ so that $a_{M\a} = a_{\a}$ as expected. 
While the ratios of $a_{\a}$ are determined by 
(\ref{c}), the normalization of $a_{\a}$ is found from (\ref{c0})
which says that the average value of $a_{\a}$ is $(h+1)/12 h^2$.
Later we give two other description of the coefficients $a_{\a}$ --- 
as certain limits and as explicit case-by-case values.   

\medskip

{\bf Conjecture.} {\em The Hirota quadratic equation 
(\ref{Hirota}---\ref{c}) coincides (up to certain rescaling
of the variables $q_k^a$) with the corresponding ADE-hierarchy of 
Kac -- Wakimoto \cite{KW}.}  

\medskip

In Section $8$ we confirm this conjecture in the cases $A_N$, $D_4$
and $E_6$. 

\medskip

\subsection*{\bf 2. The total descendent potential.}
The second goal of this paper is to generalize to the ADE-singularities
the result of \cite{GiI} that the {\em total descendent 
potential} associated to the $A_{n-1}$-singularity in the axiomatic theory 
of topological gravity is a tau-function of the $nKdV$ 
(or {\em Gelfand-Dickey}) hierarchy. 

According to E. Witten's conjecture  \cite{W} proved by M. Kontsevich 
\cite{Ko}, the following generating function for intersection indices on
the Deligne -- Mumford spaces satisfies the equation of
the KdV-hierarchy:
\footnote{Here $\psi_i$ is the $1$-st Chern class of the line bundle over 
$\M_{g,m}$ formed by the cotangent lines to the curves at the $i$-th
marked points.} 
\begin{equation} \label{WK} 
\D_{A_1} = \exp \sum_{g,m}\frac{\h^{g-1}}{m!}\int_{\M_{g,m}}
\prod_{i=1}^m (\psi_i+\sum_{k=0}^{\infty} q_k\psi_i^k ) .\end{equation}

In the axiomatic theory, the total descendent potential is, by definition,
an asymptotical function of the form
\[ \D = \exp \sum_{g\geq 0} \h^{g-1}\F^{(g)}(\q)\]
where $\F^{(g)}$ are formal functions on $\H_{+}$ {\em near the point
$\q = -1z$}. (Here $1$ is the unit element in the local algebra $H$.) This
convention called the {\em dilaton shift} is already explicitly 
present in (\ref{WK}). The formal functions $\F^{(g)}$
called the {\em genus $g$ descendent potentials} are supposed to satisfy 
certain axioms dictated by Gromov--Witten theory. The axioms (while not
entirely known) are to include the so called {\em string equation} (SE),
{\em dilaton equation} (DE), {\em topological recursion relations} (TRR or 
$3g-2$-jet property) and {\em Virasoro constraints}. 

According to \cite{GiF}, the genus-$0$ axioms SE$+$DE$+$TRR for $\F^{(0)}$
are equivalent to the following geometrical property ($\ast$) of the 
Lagrangian submanifold $\L \subset \H = T^*\H_{+}$ defined as the
graph of $d\F^{(0)}$:

{\em ($\ast$)\ $\L$ is a Lagrangian cone with the vertex at the origin and 
such that tangent spaces $L$ to $\L$ are tangent to $\L$ exactly along $zL$.}

In other words, the cone $\L$ is swept by the family 
$\tau\in H \mapsto zL_{\tau}$ of isotropic subspaces which form 
a variation of semi-infinite Hodge structures in the sense of 
S. Barannikov \cite{Ba}. According to his results, this defines
a Frobenius structure on the space of parameters $\tau$. 

In the case of ADE-singularities (and, more generally, finite reflection
groups) the Frobenius structures have been constructed by K. Saito \cite{S}.
Consider the miniversal deformation 
\[ f_{\tau}(x)=f(x) + \tau^1 \psi_1(x)+...+\tau^N \psi_N(x), \] 
where $\{ \psi_a \}$ form a weighted-homogeneous basis in the local 
algebra $H$, and $\psi_N=1$.
The tangent spaces $T_{\tau}\T$ to the parameter space $\T\simeq \CC^N$
are canonically identified with 
the algebras of functions on the critical schemes $crit (f_{\tau})$:
$\p_{\tau^a}\mapsto \p f_{\tau}/\p \tau^a \mod (\p f_{\tau}/\p x)$.   
The multiplication $\bullet $ on the tangent spaces is Frobenius with
respect to the following {\em residue metric}:
\[ (\p_{\tau_a},\p_{\tau_b})_{\tau}:= (\frac{1}{2\pi i})^3 
\oint \oint \oint \frac{\psi_a(x) \psi_b(x)\ \go }
{\frac{\p f_{\tau}}{\p x_1} \frac{\p f_{\tau}}{\p x_2} 
\frac{\p f_{\tau}}{\p x_3}} .\]
The residue metric is known to be flat and together with the Frobenius 
multiplication, the unit vectors $\p_{\tau^N}$ and the 
{\em Euler vector field}
\[ E:= \sum_{a=1}^N\ (\deg \tau^a)\ \tau^a \ \p_{\tau^a}, \ \ 
\deg \tau^{a} = 1-(m_a-1)/h, \]
forms a conformal Frobenius structure on $\T$ (see \cite{D}).
 
On the other hand, the condition ($\ast$) involves only the symplectic
structure $\Omega$ on $\H$ and the operator of multiplication by $z$ 
and thus admits the following {\em twisted loop group} of symmetries:
\[ L^{(2)}GL(H) = \{ M\in End(H)\ ((1/z))\ |\ \ M(-z)^* M(z) =1\} .\] 
According to a result from \cite{GiF}, when the Frobenius structure associated
to the cone $\L$ is semisimple, one can identify $\L$  
with the Cartesian product $\L_{A_1}\times ... \times \L_{A_1}$ of 
$N=\dim H$ copies of the cone $\L_{A_1}$ defined by 
the genus $0$ descendent potential 
$\F^{(0)}_{A_1}=\lim_{\h\to 0} \h \ln \D_{A_1}$. The identification
is provided by a certain transformation $M_{\tau}$ from (a completed
version of) $L^{(2)}GL(H)$ whose construction depends on the choice of a 
semisimple point $\tau$.

A number of results in Gromov -- Witten theory suggests that the higher
genus theory inherits the symmetry group $L^{(2)}GL(H)$ (see \cite{GiQ, GiF}).
This motivates the following construction of the total descendent 
potential of a {\em semisimple} Frobenius manifold.   

Adopt the following rules of quantization $\hat{\ }$ of 
quadratic hamiltonians. Let $\{ ..., p_{a}, ..., q_{b}, ... \}$ be a 
Darboux coordinate
system on the symplectic space $(\H, \gO)$ compatible with the polarization
$\H=\H_{+}\oplus \H_{-}$. Then
\[ (q_{a}q_{b})\hat{\ }=q_{a}q_{b}/\h, \ \ 
   (q_{a}p_{b})\hat{\ }=q_{a} \p/\p q_{b},\ \
   (p_{a}p_{b})\hat{\ }=\h \p^2/\p q_{a} \p q_{b} .\] 
This gives a {\em projective} representation of the Lie algebra 
$L^{(2)}gl(H)$ in the Fock space. The central extension is due to
\[ [\hat{F},\hat{G}]=\{ F, G \}\hat{\ } + \C (F,G)\]
where $\C$ is the cocycle satisfying  
\[ \C (p_{\a}p_{\b}, q_{\a}q_{\b} ) = \left\{ \begin{array}{l}
1 \ \text{if}\ \a\neq \b, \\ 2\ \text{if}\ \a=\b \end{array} \right. \]
and equal $0$ for any other pair of quadratic Darboux monomials.   

Introduce the {\em total descendent potential} as an asymptotical function: 
\[ \D : = C(\tau)\ \hat{M}_{\tau}\ [\D_{A_1}\otimes...\otimes \D_{A_1}],\]
where $\hat{M}:= \exp (\ln M)\hat{\ }$, and $C(\tau)$ is a normalizing 
constant possibly needed to keep the R.H.S. independent of the choice
of a semisimple point $\tau$. This definition has been tested in 
\cite{GiQ, GiI} and is known to agree with the TRR, SE, DE and the Virasoro
constraints. Here is a more explicit description of $M_{\tau}$ and
$C(\tau)$ in the form applicable to Frobenius manifolds of simple 
singularities.
  
 Consider the complex oscillating integral
\[  \J_{\Gb}(\tau) = (- 2\pi  z)^{-3/2}\int_{\Gb} e^{f_{\tau}(x)/z} \go .\]
Here $\Gb $ is a non-compact cycle from the relative homology group
\[ \lim_{C\to \infty}\ 
H_m (\CC^m, \{ x: \operatorname{Re} (f_{\tau}/z) \leq - C \})\simeq \ZZ^N. \]
We will use the notation $\p_1,...,\p_N$ for partial derivative with 
respect to a flat (and weighted - homogeneous)
coordinate system $(t^1,...,t^N)$ of the residue metric. We treat the 
derivatives $z\p_{a} \J_{\Beta}$ as components of a covector field
$z \sum \p_a \J_{\Beta} dt^a \in T^*\T$ which can be identified
with a vector field via the residue metric and --- via its Levi-Civita
connection --- with an $H$-valued function $J_{\Beta}(z,\tau)$. 
According to 
K. Saito's theory these functions satisfy in flat coordinates 
the differential equations
\begin{equation} \label{qde} z \p_a J = (\p_a\bullet) J \end{equation} 
together with the homogeneity condition:
\begin{equation} \label{euler} (z\p_z -\mu +z^{-1}E\bullet) J=0 \end{equation}
where $\mu=-\mu^*$ is the diagonal operator with the eigenvalues  
$ 1/2-m_a/h$. 
The latter equation yields an {\em isomonodromic} family of 
connection operators $\nabla_{\tau}=\p_z-\mu/z+(E\bullet)/z^2$
regular at $z=\infty$ and turning into $\p_z - \mu /z$ at $\tau=0$.

According to \cite{D2}, there exists a (unique in the ADE-case)
gauge transformation of the form $S_{\tau}(z)=\1+S_1(\tau)z^{-1}+S_2(\tau)
z^{-2}+...$ ({\em i.e.} near $z=\infty$) 
conjugating $\nabla_{\tau}$ to $\nabla_0$ and such that 
$S^*_{\tau}(-z)S_{\tau}(z)=\1$. It satisfies the homogeneity 
condition $(z\p_z +L_{E}) S_{\tau} = [\mu , S_{\tau}]$.

On the other hand, let $\tau$ be semisimple. Then the 
functions $f_{\tau}$ have $N$ non-degenerate critical points $x^{(a)}(\tau)$
with the critical values $u^a(\tau)$ and the Hessians $\Delta_a(\tau)$.
The local coordinate system $\{ u^a \}$ (called {\em canonical})
diagonalizes the product 
$\bullet$ and the residue metric: 
\[ \p/\p u^a \bullet \p/\p u^b =\d_{ab} \p/\p u^b,\ \ 
(\p /\p u^a, \p/\p u^b)_{\tau} = \d_{ab} \Delta_a^{-1} \p/\p u^b .\]    
Define an orthonormal coordinate system 
\[ \Psi(\tau): \CC^N \to T_{\tau}\T = H, \ \ \Psi (q^1,...,q^N) 
= \sum q^a \sqrt{\Delta_a} \p/\p u^a, \]
and put $U_{\tau}=\diag [u^1(\tau),...,u^N(\tau)]$.
Stationary phase asymptotics of the oscillating integrals 
$J_{\Beta_a}, a=1,...,N$, near the corresponding critical points $x^{(a)}$
yield a fundamental solutions to the system (\ref{qde}),(\ref{euler}) 
in the form
\[ \Psi(\tau) R_{\tau}(z) e^{U_{\tau}/z}, \ \ 
R_{\tau}=\1+R_1(\tau)z+R_2(\tau)z^2+...,\ \ R^t_{\tau}(-z)R_{\tau}(z)=\1 .\]
The matrix series $R_{\tau}$ satisfies the homogeneity condition 
$(z\p_z+L_{E})R_{\tau} =0$ and, according to \cite{GiQ}, an
asymptotical solution with this property is unique up to reordering
or reversing the basis vectors in $\CC^N$. 

Define
\[ c(\tau) := \frac{1}{2}\int^{\tau}\sum_a R_1^{aa}(\tau)du^a \]
as the local potential of the $1$-form $\sum R_1^{aa}du^a/2$ 
(which is known to be closed \cite{D}).  

In the above notations, the total descendent potential of the
ADE-singularity assumes the form
\begin{equation} \label{D} \D = e^{c(\tau)}\ \hat{S}_{\tau}^{-1}\ 
\Psi(\tau)\ \hat{R}_{\tau}\ e^{\hat{\frac{U_{\tau}}{z}}} \  
\D_{A_1}^{\otimes N} .
\end{equation}
The R.H.S. is known to be independent of $\tau$ (see \cite{GiQ})
and defines $\D$ (up to a constant factor) as an asymptotical function
of $\q=\q_0+\q_1z+\q_2z^2+...$ in the formal neighborhood
of $\q = \tau-z$ with semisimple $\tau$.

Our main result is the following theorem.

\medskip

{\bf Theorem $1$.} {\em The total descendent potential (\ref{D}) of a
simple singularity satisfies the corresponding Hirota quadratic 
equation (\ref{Hirota}--\ref{c}).}

\medskip
  
In Section $4$, we discuss Hirota quadratic equations of the 
$KdV$-hierarchy. The plan for the proof of Theorem $1$ is to reduce 
the Hirota quadratic equations for $\D$ to those for $\D_{A_1}$ 
by conjugating the vertex operators
in (\ref{Hirota},\ref{L-nod}) past the quantized symplectic transformations
from (\ref{D}). In Section $5$, we describe the results of such 
conjugations by quoting corresponding theorems from \cite{GiI}. The
residue in (\ref{Hirota}) is computed in Section $6$ and is compared with
(\ref{L-nod}) in Section $7$. The case-by-case tables
for the coefficients $a_{\a}$ are presented in Section $8$.
 A key to all our computations is the phase form and its properties 
discussed in next section.

\medskip

\subsection*{\bf 3. The phase forms and the root systems.} Consider a flat
family of cycles $\phi \in H_2(f_{\tau}^{-1}(\gl))$ in the non-singular 
Milnor fibers and define the 
{\em period vector} $I^{(0)}_{\phi} (\gl,\tau) \in H$ by
\begin{equation} \label{period} 
(I^{(0)}_{\phi}(\gl,\tau), \p_a) :=\p_a\ \frac{(-1)}{2\pi} 
\int_{\phi \subset f_{\tau}^{-1}(\gl)} \frac{\go}{df_{\tau}}. \end{equation} 
It is a multiple-valued vector function on the complement to the discriminant
which turns into $I^{(0)}_{\phi}(\gl)$ from Section $1$ at $\tau =0$.

The {\em phase form} $\tilde{\W}_{\a,\b}$ (defined in \cite{GiI}, Section $7$)
is given by the formula
\[ \tilde{\W}_{\a,\b}(\gl,\tau) := \sum_{i=1}^N 
(I^{(0)}_{\a}(\gl,\tau) \bullet I^{(0)}_{\b}(\gl,\tau), \p_{\tau^a}) 
\ d\tau^a. \]
It is a multiple-valued $1$-form on the complement to the discriminant
and depends bilinearly on the cycles $\a,\b$ (to be chosen in 
$(f_0^{-1}(1)$ and transported to $f_{\tau}^{-1}(\gl)$).
According to \cite{GiI}, the phase forms have the following properties.
\begin{enumerate}
\item $d\tilde{\W}_{\a,\b}=0$. 
\item $L_{\p_{\gl}+\p_N} \tilde{\W}_{\a,\b} =0$, {\em i.e.} $\tilde{\W}$
is determined by the restriction 
\[ \W_{\a,\b}(\tau):=\tilde{\W}_{\a,\b}(0,\tau), \ \  
\tilde{\W}_{\a,\b}(\gl,\tau)=\W_{\a,\b}(\tau-\gl \1) .\]
\item $L_E \W_{\a,\b} = 0$.
\item Near a generic point of the discriminant $\Delta \subset \T$
the form $\W_{\a,\b}$ becomes single-valued on the double cover and
has a pole of order $\leq 1$ on $\D$ (since $I^{(0)}_{\a}$ have a pole
of order $\leq 1/2$).
\item $\oint_{\d_{\c}} \W_{\a,\b} = -2\pi i \lan \a, \c\ran \lan \b,\c\ran,$ 
where $\c$ is the cycle vanishing over a generic point of the discriminant,
and $\d_{\c}$ is a small loop going {\em twice} 
(in the positive direction defined by complex orientations) 
around the discriminant near this point.
\end{enumerate}

\medskip

{\bf Proposition 1.} 
\[ \W_{\a\b} = -\frac{1}{2}\sum_{\c\in A} \lan\a,\c\ran\lan\b,\c\ran
\frac{d\lan\c,x\ran}{\lan\c,x\ran} .\]  

\medskip

{\em Proof.} The phase form $\W_{\a,\b}$ becomes single-valued 
on the {\em Chevalley cover} representing $\T$ as the quotient
of $\CC^N = H^2(f_0^{-1}(1),\CC)$ by the monodromy group.
The properties ($1$) and ($4$) show that it has at most 
logarithmic pole on the mirrors $\lan \c,x\ran =0$. The property
($5$) controls the residues on the mirrors. The difference
of the L.H.S. and the R.H.S. has to be a holomorphic $1$-form,
homogeneous of degree $0$ by to the property ($3$), and therefore
vanishes identically. $\square$

\medskip

{\bf Corollary 1.} $i_E \W_{\a,\b} = -\lan \a,\b\ran $.

\medskip

Indeed, the Euler vector field becomes $h^{-1}\sum x_a \p_{x_a}$ on the
Chevalley cover, so that the equality follows from 
$\sum_{\c\in A} \lan \a,\c\ran \lan \b,\c\ran = 2h \lan\a,\b\ran$.
This is one more general property of phase forms established in \cite{GiI}.

\medskip

{\bf Corollary 2.} {\em The phase form $\W_{\b}$ of Section $1$
coincides with $\W_{\b,\b}$.}

\medskip

{\em Remark.} The inverse to the Chevalley quotient map is given
by the period map 
\[ \tau \mapsto [\go /df_{\tau}] \in H^2(f_{\tau}^{-1}(0),\CC) \leadsto
H^2(f_0^{-1}(1),\CC) \simeq \CC^N.\]
The periods $I_{\a}^{(0)}$ are defined via the differential of the 
inverse Chevalley map and therefore represent parallel translations
of the cycles $\a$ considered as {\em co}vectors in $\CC^N$. 
The value of phase form $\W_{\a,\b}$, which is also a covector, 
is constructed as the Frobenius product $\a \bullet \b$ of {\em co}vectors 
(defined by the isomorphisms $T_{\tau}\T \simeq T^*_{\tau}\T$ based
on the residue metric). Thus the formula
\[ \a \bullet \b := 
\frac{1}{2} \sum_{\c\in A} \frac{\lan\a,\c\ran\lan\b,\c\ran}
{\lan \c,x\ran}\ \c \]
defines on $(\CC^N)^*$ a family of commutative associative 
multiplications depending on the parameter $x$.
It would be interesting to find a representation-theoretic
interpretation of this structure defined entirely in terms of 
the root system $A$.   

\medskip

We prove several further properties of phase forms needed in our 
computations.

\medskip

{\bf Proposition 2.} {\em In the case of $ADE$-singularities, suppose
that $\b$ has integer intersection indices with all $\a\in A$ and 
is invariant under the monodromy around a discriminant-avoiding loop $\c$. 
Then $\oint_{\c} \W_{\b,\b} \in 2\pi i\ \ZZ$.}

\medskip

{\em Proof.} This is Proposition $1$ from Section $7$ of \cite{GiI}.
$\square$

It would be interesting to find out if the property remains valid 
for non-simple singularities.

\medskip

We will see in Section $7$ that the coefficients $a_{\a}$ introduced
in Section $1$ can be equivalently defined via the following limits
$b_{\a}$. Start with choosing $(\tau_1,...,\tau_N)=-\1=(0,...,0,-1)$ 
in the role of the base point in $\T$ and identify $A$ with the set 
of vanishing cycles in $H_2(f_{-\1}^{-1}(0)) = H_2(f_0^{-1}(1))$.
Let us also fix $\tau\in \T$ such that $f_{\tau}$ is a Morse function,
and let $u$ will be one of the critical values of $f_{\tau}$ so that
$\tau-u\1 \in \Delta$. We may assume
that $\tau -(u+1)\1 \notin \Delta$ and that the straight segment connecting 
$\tau -(u+1)\1$ with $\tau-u\1$ does not intersect $\Delta$. 
For each $\a\in A$,
pick a discriminant-avoiding path $\c_{\a}$ connecting $-1$ with 
$\tau-(u+1)\1$ and further with $\tau-u\1$ along the straight segment and
such that $\a$ becomes the vanishing cycle when transported along $\c_{\a}$ 
from $\1$ to $\tau-u\1$. Assuming that integration of the phase form
is performed along this path we put
\begin{equation} \label{b}
 b_{\a} :=\lim_{\ge\to 0}\ \exp \left\{\ -\
\int_{-\1}^{\tau-(u+\ge)\1} \W_{\a,\a}\ - 
\ \int_{-1}^{-\ge}\frac{2\ dt}{t} \ \right\} % = \lim_{\ge\to 0}
%\exp \left\{ \int^{-\ge\1}_{\tau-(u+\ge)\1} \W_{\a,\a} \right\} .
\end{equation}

{\bf Proposition 3.} 
{\em The limit exists and does not depend on the choice of the
path of integration provided that the path terminates at a generic
point of the discriminant and that the cycle $\a$ transported along the
path vanishes over this point.}

\medskip

{\em Proof.} We may assume that $u=u^1$ is the first of the 
canonical coordinates $U=(u^1,...,u^N)$, and therefore $u^1=0$ 
is the local equation of the discriminant branch. Since $\a$ is vanishing
at the end of the path, the period vector $I^{(0)}_{\a}$ has the following
expansion (here $\1_i$ stand for the standard basis vectors in $\CC^N$): 
\[ \Psi^{-1}_{\tau}I^{(0)}_{\a}(\gl,\tau)=\frac{\pm 2}{\sqrt{2(\gl-u^1)}}
\left( \1_1 + (\gl-u^1) \sum a^i(U) \1_i+o (\gl-u^1) \right) .\]
Since $\Psi (\1_i) = \sqrt{\Delta_i} \p/\p u^i$, we have
$(\1_i\bullet \1_j, \p/\p u^k) = \d_{ij}\d_{ik}$, and therefore 
\[ \W_{\a,\a} = \sum (I^{(0)}_{\a}\bullet I^{(0)}_{\a}, \p/\p u^k) du^k
\ \left|_{\gl=0} \right. = \frac{2 du^1}{-u^1} + 4 a^1(U) du_1 + 
{\mathcal O} (-u_1).\]
We see that the integral $\int_{u_1=-1}^{u_1=0} \W_{\a,\a}$ diverges 
the same way as $-\int_{-1}^0 2dt/t$ so that the difference converges.
This proves the existence of the limit. Removing this singular term
we find that the integral $\int [4 a_1(U) du^1 + {\mathcal O} (-u_1)]$ 
vanishes along any path inside the discriminant branch $u_1=0$. 
This shows that
the limit $b_{\a}$ is locally constant as a function 
of the path's endpoint on the discriminant, and therefore --- globally constant
due to the irreducibility of the discriminant. Finally, precomposing a path
with a discriminant-avoiding loop $\c$ with trivial monodromy of the
cycle $\a$ does not change $b_{\a}$ thanks to Proposition $2$. $\square$

\medskip

{\bf Corollary.} {\em $a_{\a}/a_{\b} = b_{\a}/b_{\b}$ for all $\a,\b\in A$.}

\medskip

{\bf Proposition $4$.} {\em Let $\d_{\ge}$ be a small loop of radius $\ge$ 
around the discriminant near a generic point $\tau-u\1$, and let 
$\lan \a,\b\ran =\pm 1$, where $\b$ is the
cycle vanishing at this point. Then $\lim_{\ge\to 0} \oint_{\d_{\ge}} 
\W_{\a,\a} = -\pi i$.}

\medskip

{\em Proof.} We have $\a=\pm \b/2+\a'$ where $\a'$ is invariant under the 
monodromy around $\d_{\ge}$. Expanding $I^{(0)}_{\a} = 
I^{(0)}_{\a'}\pm I^{(0)}_{\b/2}$
near $\gl=u$ as in the proof of Proposition $3$ we find 
\[ \oint_{\d_{\ge}}\W_{\a,\a} = \oint_{\d_{\ge}}\frac{du}{-2u} + 
\oint_{\d_{\ge}}{\mathcal O} (\sqrt{u})\ du = -\pi i + 
{\mathcal O} (\sqrt{\ge}) \to -\pi i. \hspace{1cm} \square \] 

In fact this property has been already used in \cite{GiI}. 

\medskip

\subsection*{\bf 4. Two forms of the KdV-hierarchy.} 
Consider the miniversal deformation
of the $A_1$-singularity in the form $f_u(x):=(x_1^2+x_2^2+x_3^2)/2+u$. 
The vanishing cycle $\a$ can be identified with the real
sphere $(x_1^2+x_2^2+x_3^2) =2(\gl-u) $. The period 
\[ \int_{\a} \frac{dx_1\w dx_2 \w dx_3}{df_u} = \frac{d}{d\gl} \frac{4}{3}
\pi (2(\gl-u))^{3/2} = 4\pi \sqrt{2(\gl-u)} .\]  
Since $(1,1)=\Res dx_1\w dx_2\w dx_3 /x_1x_2x_3 = 1$, we have
$ I_{\a}^{(-1)}(\gl,u) = 2 \sqrt{2(\gl-u)}$, and more generally,
$ I^{(k)}_{\pm \a}(\gl,u) = \pm 2 (d/d\gl)^k (2(\gl-u))^{-1/2}, \ k\in \ZZ$.
The Coxeter transformation swaps $\a$ and $-\a$ and so 
$a_{\a} = a_{-\a} =  (h+1)/12h^2 = 1/16$. The equation 
(\ref{Hirota},\ref{L-nod}) in this example assumes the form 
\begin{equation} \label{SL_2}   
\Res_{\gl =\infty} \frac{d\gl}{\gl} \left[ \sum_{\pm} 
\ ^{A_1}\Gamma^{\pm\a}(\gl)\otimes \ ^{A_1}\Gamma^{\mp \a}(\gl) \right]
(\Phi \otimes \Phi) = 16\ (l + \frac{1}{8})\ (\Phi \otimes \Phi) ,
\end{equation}
where
\begin{equation} \label{l} 
l = \sum_{k\geq 0} \frac{2k+1}{2} (q_k\otimes 1-1\otimes q_k)\ 
(\p_{q_k}\otimes 1 - 1\otimes \p_{q_k}) .\end{equation}
Here we use the notation $\ ^{A_1}\Gamma^{\phi}(\gl)$ to single out the
vertex operators $\Gamma^{\phi}(\gl)$ of the $A_1$-singularity. 
In order to identify the condition (\ref{SL_2}) for $\Phi$ 
with the KdV hierarchy in \cite{K,KW} 
corresponding to the root system $A_1$, we denote 
$\sqrt{2\gl}$ by $\zeta$, rescale the variables by $q_k = (2k+1)!! t_{2k+1}$ 
and put $x_m =(t'_m+t''_m)/2$, $y_m=(t'_m-t''_m)/2$ where $m=1,3,5,...$.
In this notation $l= \sum m y_m \p_{y_m}$, and (\ref{SL_2},\ref{l}) 
becomes
\[ \left[  \Res \frac{d\zeta}{\zeta}\ e^{4\sum  
\zeta^m \frac{y_m}{\sqrt{\h}}} \ e^{-2\sum \frac{\zeta^{-m}}{m}
\sqrt{\h}\p_{y_m}}
-1-8 \sum m y_m\p_{y_m} \right]\ \Phi(\x+\y)\ \Phi(\x-\y) = 0 .\]     
This coincides with the equation $(14.13.1)$ in \cite{K} characterizing 
tau-functions $\Phi $ of the KdV hierarchy. 

Another form of the Hirota quadratic equation for $\Phi$ is based
on the representation of the KdV-hierarchy as the $\mod 2$-reduction of 
the KP-hierarchy (see \cite{K}, Section $14.11$). 
It can be rephrased (see \cite{GiI}) as the condition 
\begin{equation} \label{KdV} 
\left[ \sum_{\pm}\ ^{A_1}\Gamma^{\pm\a/2}(\gl)\otimes 
\ ^{A_1}\Gamma^{\mp\a/2}(\gl)\frac{d\gl}{\pm \sqrt{\gl}}
\right] \ (\Phi\otimes \Phi)
\ \text{\em has no pole in $\gl$}.\end{equation}
Indeed, in the previous notations this can be rewritten as the
property 
\[ e^{2\sum  
\zeta^m \frac{y_m}{\sqrt{\h}}} \ e^{-\sum \frac{\zeta^{-m}}{m}
\sqrt{\h}\p_{y_m}} \ \Phi(\x+\y) \Phi(\x-\y) \ \text{\em contains no 
$\zeta^{-m}$ for odd $m>0$}.\]
This coincides with the $\mod 2$-reduction of the KP-hierarchy of
the Hirota equation $(14.11.5)$ in \cite{K}. 
According to a result from \cite{K}, Section $14.13$, this condition is
actually equivalent to (\ref{SL_2}). 

In Section $6$ we will use the fact that (according to Kontsevich's theorem)
the function $\Phi = \D_{A_1}$ satisfies both forms (\ref{SL_2}) and 
(\ref{KdV}) of the KdV-hierarchy.

\medskip
    
\subsection*{\bf 5. Symplectic transformations of vertex operators.}
Generalizing the construction of Section $1$, introduce the vertex
operator $\Gamma^{\phi}_{\tau}(\gl)$ corresponding to the vector 
$\f \in H[[z,z^{-1}]]$ of the form  
\[ \f^{\phi}_{\tau}(\gl):=\sum_{k\in \ZZ} I^{(k)}_{\phi}(\gl,\tau) (-z)^k .\] 
Here $I^{(0)}_{\phi}$ is the period vector introduced in Section $3$,
and $I^{(k)}_{\phi}:= d^k I^{(0)}_{\phi}/d\gl^k$ as before. For $k<0$
the integration constants are taken ``equal $0$'' so that 
$I^{(k)}_{\phi}$ satisfy the homogeneity conditions:
\[  (\gl \p_{\gl} + L_{E})\ I^{(k)}_{\phi}(\gl,\tau) = (\mu-\frac{1}{2}-k)\ 
I^{(k)}_{\phi}(\gl,\tau) .\]
In particular $\Gamma_0^{\phi}$ coincides with the
vertex operator $\Gamma^{\phi}$ from Section $1$.
We state below several results about behavior of the vertex operators
under conjugation by some symplectic transformations and refer to 
Sections $5,6,7$ of \cite{GiI} for the proofs.

\medskip

{\bf Theorem A} (see Proposition $2$ in \cite{GiI}). 
\[ \hat{S}_{\tau}\ \Gamma^{\phi}_0(\gl) \ \hat{S}_{\tau}^{-1}\ = \ 
\exp\left\{ \frac{1}{2}\int_{-\gl\1}^{\tau-\gl\1} \W_{\phi,\phi} \right\}\ \  
\Gamma^{\phi}_{\tau}(\gl) \]

\medskip

We have to stress here that in order to compare the vertex operators 
$\Gamma_{0}^{\phi}(\gl)$ and $\Gamma_{\tau}^{\phi}(\gl)$ one needs
to transport the cycle $\phi$ from $f_{0}^{-1}(\gl)$ to $f_{\tau}^{-1}(\gl)$
along a path in $\T$ connecting $-\gl\1=(0,...,-\gl)$
with $\tau-\gl \1=(\tau_1,...,\tau_N-\gl)$ and avoiding the discriminant
$\Delta$ corresponding to singular levels $f_{\tau}^{-1}(0)$.
It is assumed in the formulation of the theorem that the integral 
of the phase form is taken along this very path. Similar conventions apply to 
other formulas of this and
following sections involving integration of phase forms. 

Now let the cycle $\phi \in H_2(f_{\tau}^{-1}(\gl)$ be written as the
sum $\phi = \lan \phi, \b\ran \b/2 + \phi'$ where $\lan \phi',\b\ran =0$.
Here $\b$ is the cycle
{\em vanishing} at a non-degenerate critical point of the function $f_{\tau}$
with the critical value $u$ and transported to $f_{\tau}^{-1}(\gl)$
along a discriminant-avoiding path connecting $\tau-\gl\1$ and $\tau-u\1$. 

\medskip

{\bf Theorem B} (see Proposition $4$ in \cite{GiI}).
\[ \Gamma_{\tau}^{\phi}(\gl) = 
\exp\left\{
\frac{\lan \phi,\b\ran}{2}\int_{\tau-\gl\1}^{\tau-u\1} \W_{\b,\phi'}\right\}\ 
\ \Gamma_{\tau}^{\phi'}(\gl)\ 
\Gamma_{\tau}^{\frac{\lan \phi,\b\ran}{2} \b} (\gl)\]

\medskip

The integral here is taken along the path terminating
on the discriminant where the phase form is singular. However the singularity
is proportional to $(\gl-u)^{-1/2}$ and is therefore integrable.

Let us recall that the columns of the matrix $R_{\tau}$ in the asymptotical
expansion $\Psi(\tau)R_{\tau}(z)\exp (U/z)$ correspond to non-degenerate
critical points of the Morse function $f_{\tau}$ with the critical values
$u^i(\tau)$. Let $\b_i$ be the cycle vanishing over $u^i$. 

\medskip

{\bf Theorem C} (see Proposition $3$ in \cite{GiI}).
\[ (\Psi(\tau) \hat{R}_{\tau})^{-1}\ \Gamma_{\tau}^{c\b_i}(\gl) \ 
(\Psi(\tau) \hat{R}_{\tau}) = e^{c^2 W_i/2} \ [ \cdots \1\otimes 
(\ ^{A_1}\Gamma^{c\b}_{u^i}(\gl))^{(i)}\otimes \1 \cdots ],\]
{\em where 
\[ W_i := \int_{\tau-\gl\1}^{\tau-u^i\1} 
\left( \W_{\b_i,\b_i}-\frac{2\ dt^N}{\tau^N-u^i-t^N} \right),\]
and $\ ^{A_1}\Gamma_u^{c\b}(\gl)$ is the vertex operator of the 
$A_1$-singularity with the miniversal deformation 
$\frac{x_1^2}{2}+\frac{x_2^2}{2}+\frac{x_3^2}{2} +u$
corresponding to the $c$-multiple of the vanishing cycle.}  

\medskip

The behavior of $I^{(0)}_{\b_i}$ near $\gl=u^i$ is described by the 
asymptotics
\[  \Psi^{-1}(\tau)\ I^{(0)}_{\b_i}(\gl,\tau) = \frac{2}{\sqrt{2(\gl-u^i)}} 
\ \left( \1_i + ... \right) \]
where $\1_i$ is the $i$-th basis vector in $\CC^N$, and the dots mean higher
order powers of $\gl-u^i$. Respectively, the vertex operator of the 
$A_1$-singularity is more explicitly defined by the series 
$\f \in \CC [[z,z^{-1}]]$ of the form
\[ \f = \sum_{k\in \ZZ} \frac{d^k}{d\gl^k} \frac{2c}{\sqrt{2(\gl-u)}}\ (-z)^k,
 \]
where the branch of the square root should be the same as in the above 
asymptotics. The subscript $(i)$ indicates the position of the vertex operator
in the tensor product operator acting on the Fock space of functions
of $(\q^{(1)},...,\q^{(N)}) = \Psi^{-1}(\tau) \q$. 
The integrand in the formula for $W_i$ considered as a $1$-form in the 
space with coordinates $(t^1,...,t^N)$ identical to parameters of the 
miniversal deformation, while the notation $\tau=(\tau^1, ..., \tau^N)$
is reserved for expressing the limits of integration. The phase form $\W$
has a non-integrable singularity at $t=\tau-u^i\1$ which happens to 
cancel out with that of the subtracted term so that the difference is
integrable.

Finally, the following result is the special case of Theorem $A$ corresponding
to the $A_1$-singularity. 
  
\medskip

{\bf Theorem D} (see Proposition $3$ in \cite{GiI}).
\[ e^{-(u/z)\hat{\ }}\ ^{A_1}\Gamma_u^{c\b}(\gl)\ e^{(u/z)\hat{\ }} = 
\exp \left\{ -\frac{c^2}{2} \int_{\gl-u^i}^{\gl} \frac{2\ dt}{t} \right\}
\ ^{A_1}\Gamma_0^{c\b}(\gl)\]  

\medskip

In fact this result can be obtained more directly using Taylor's formula.
Indeed, for any analytic function $I^{(0)}$ we have
\[ e^{-u/z} \left[ \sum_{k\in \ZZ} I^{(k)}(\gl) (-z)^k \right] e^{u/z} =
\sum_{k\in \ZZ} I^{(k)}(\gl+u) (-z)^k \]
provided that $|u|$ does not exceed the convergence radius of $I^{(0)}$
at $\gl$. Thus the transformation in the theorem effectively 
consists in the translation $\sqrt{\gl-u} \leadsto \sqrt{\gl}$ 
along an origin-avoiding path. The
integral in the exponent should be taken along this path.       
  
\medskip

\subsection*{\bf 6. The residue sum.} In this section, we compute
the residue sum 
\begin{equation} \label{res} 
\Res_{\gl=\infty} \frac{d\gl}{\gl} \left[ \sum_{\a\in A} 
b_{\a} \Gamma_0^{\a}(\gl)\otimes \Gamma_0^{-\a}(\gl) \right]\ 
\D^{\otimes 2} \end{equation}
assuming that the coefficients $b_{\a}$ are defined as in Proposition $3$. 

Introduce the {\em total ancestor potential}
\[ \A_{\tau}:=\ \hat{S}_{\tau}\ \D\ =\ e^{c(\tau)}\ \Psi(\tau)\ 
\hat{R}_{\tau}\ e^{(U_{\tau}/z)\hat{\ }}\ \D_{A_1}^{\otimes N} .\]
Applying Theorem A of the previous section we find that (\ref{res})
can be rewritten as
\begin{equation} \label{res1} 
\Res_{\gl=\infty} \gl d\gl \left[ \sum_{\a\in A} 
c_{\a} \Gamma_{\tau}^{\a}(\gl)\otimes \Gamma_{\tau}^{-\a}(\gl) \right]\ 
\A_{\tau}^{\otimes 2}, \end{equation}
where 
\[ c_{\a} =\lim_{\ge\to 0}\  \exp \left\{\ -\ 
\int_{\tau-\gl\1}^{\tau-(u+\ge)\1} \W_{\a,\a}\ - 
\ \int_{-1}^{-\ge}\frac{2\ dt}{t} \ \right\}. \]
assuming that $\a \in H_2(f_{\tau}^{-1}(\gl))$ vanishes at $\gl=u$ when
transported along the path of integration of the phase form. 
Note that the factor $\gl^{-1}d\gl$ in (\ref{res}) is replaced by $\gl d\gl$
in (\ref{res1}) due to Corollary $1$ from Section $3$ which shows that 
\[ \exp \left\{ -\int_{-\1}^{-\gl\1} \W_{\a,\a} \right\} = 
\exp \left\{ \lan \a,\a\ran \int_{-1}^{-\gl} \frac{dt}{t} \right\} = \gl^2.\]  

The ancestor potential $\A_{\tau} = \exp \sum \h^{(g-1)} \F_{\tau}^{(g)}$ 
is a {\em tame} asymptotical function in the following sense: 
$\F^{(g)}_{\tau}$ considered as a formal function of $t_k^a=q_k^a + 
\d_{k1}\d_{aN}$ satisfy 
\[  \frac{\p^r \F^{(g)}_{\tau}}{\p t_{k_1}^{a_1} ... \p t_{k_r}^{a_r}} 
|_{\t =0} = 0 \ \ \text{whenever} \ \ k_1+...+k_r > 3g-3+r .\]
This follows from the analogous property of $\D_{A_1}^{\otimes N}$,
from the invariance of $\D_{A_1}$ under the string flow $\exp (u/z)\hat{\ }$ 
and from the ``upper-triangular'' property of $\hat{R}_{\tau}$.
We refer to Proposition $5$ in \cite{GiI} for the proof.  
It is also shown in Section $8$ of \cite{GiI} that for tame asymptotical
functions $\Phi $ the vertex operator expressions 
$\Gamma_{\tau}^{\phi}(\gl) \otimes \Gamma_{\tau}^{-\phi} (\gl) \ 
\Phi^{\otimes 2}$ can be considered not only as series expansions 
in fractional powers of $\gl$ near $\gl = \infty$, but also as 
multiple-valued analytical functions defined over the entire range of $\gl$
and ramified only on the discriminant.
Moreover, the sum in (\ref{res1}) is manifestly invariant under the entire
monodromy group ($=$ the ADE-reflection group). Therefore the sum is actually
a single-valued differential $1$-form on the complement to $\D$. 
Thus the residue (\ref{res1}) at $\gl =\infty$ coincides with the sum
of residues at the critical values $\gl=u_i$ of the function $f_{\tau}$.
Our next goal is to take $u=u_i$ and compute the residue. 
    
In a neighborhood of $\gl=u$, the monodromy group reduces to $\ZZ_2$
generated by the reflection $\sigma$ in the hyperplane orthogonal to 
two vanishing cycles which we denote $\pm \b$. 

\medskip

First, consider the summand in (\ref{res1}) corresponding to a 
$\sigma$-invariant cycle $\a\in A$. The period vectors 
$I^{(k)}_{\a}(\gl,\tau)$ are therefore single-valued analytic functions
near $\gl=u$. In particular, $\ln c_{\a}$, which differs from a constant
by $\int_{\tau-(u+1)\1}^{\tau-\gl\1}\W_{\a,\a}$, is analytic
too. We conclude that 
$\gl\ c_{\a}\ \Gamma_{\tau}^{\a}(\gl) \otimes \Gamma_{\tau}^{-\a}(\gl)\  
\A^{\otimes 2}$ has no pole at $\gl=u$. 

\medskip 

Next, consider a pair of cycles $\a_{\pm}\in A$ transposed by $\sigma$
and having intersection indices $\pm 1$ with $\b$. We have
$\a_{\pm}=\a'\pm \b/2$ where $\sigma \a'=\a'$.
We use Theorem B to replace $\Gamma_{\tau}^{\a_{\pm}}$ with 
$\Gamma_{\tau}^{\a'} \Gamma_{\tau}^{\pm \b/2}$ and then commute 
$\Gamma_{\tau}^{\pm \b/2}$ across $\Psi \hat{R} \exp (U/z)\hat{\ }$ using 
Theorems C and D. The terms from (\ref{res1}) corresponding to
$\a=\a_{\pm}$ turn into
\begin{align} \label{res2} \gl\ d\gl \ 
\left[ \Gamma_{\tau}^{\a'}(\gl)\otimes \Gamma_{\tau}^{-\a'}(\gl) \right]\ 
\left( \Psi(\tau) \hat{R}_{\tau} e^{(U_{\tau}/z)\hat{\ }}\ \otimes \ 
\Psi(\tau) \hat{R}_{\tau} e^{(U_{\tau}/z)\hat{\ }} \right) \\ \notag
\times \ \left[  ... \1 \otimes \left( \sum_{\pm} d_{\pm}
\ ^{A_1} \Gamma_{0}^{\pm\b/2}(\gl) \otimes \ ^{A_1}\Gamma_{0}^{\mp\b/2}(\gl) 
\right)^{(i)} \otimes \1  ... \right] \ 
\left( \D_{A_1}^{\otimes N} \otimes \D_{A_1}^{\otimes N}\right) .
\end{align}
The coefficients $d_{\pm}$ here are   
\begin{align} \notag  d_{\pm} = \lim_{\ge\to 0} \exp \left\{ 
- \oint_{\tau-\gl\1}^{\tau-(u+\ge)\1}
\W_{\a_{\pm},\a_{\pm}} - \int_{-1}^{-\ge}\frac{2dt}{t}
\pm \int_{\tau-\gl\1}^{\tau-u\1}\W_{\a',\b} + \right. \\ \label{a}
\left. \int_{\tau-\gl\1}^{\tau-(u+\ge)\1}\W_{\b/2,\b/2}+
\int_{u-\gl}^{-\ge}\frac{dt}{2t} - 
\int_{u-\gl}^{-\gl}\frac{dt}{2t} \right\}.
\end{align}
We have to emphasize that all integrals here except the
first one are taken along a short path near $\gl=u$ making $\b$ vanish
while in the first integral this path is precomposed with a loop 
transforming $\a_{\pm}$ to $\b$. 

Let us take $\gl=u+1$ for the base point for such a loop $\c_{\pm}$ and
rearrange the first integral as
\[ -\int_{\c_{\pm}} \W_{\a_{\pm},\a_{\pm}}
-\int_{\tau-(u+1)\1}^{\tau-(u+\ge)\1} \W_{\b,\b}  
+ \int_{\tau-(u+1)\1}^{\tau-\gl\1} \W_{\a_{\pm},\a_{\pm}} .\]  
Combining this with $\W_{\a_{\pm},\a_{\pm}}=\W_{\a',\a'}\pm \W_{\a',\b}+
\W_{\b/2,\b/2}$ we can rewrite the exponent in (\ref{a}) as
\begin{align} \label{1} 
- \int_{\c_{\pm}} \W_{\a_{\pm},\a_{\pm}} 
- \int_{ \tau-(u+1)\1}^{\tau-(u+\ge)\1}\W_{\b,\b} 
- \int_{-1}^{-\ge} \frac{2dt}{t} 
\pm \int_{\tau-(u+1)\1}^{\tau-(u+\ge)\1} \W_{\a',\b} \\
\label{2}
\int_{\tau-(u+1)\1}^{\tau-\gl\1} \W_{\a',\a'} + 
\int_{\tau-(u+1)\1}^{\tau-(u+\ge)\1} \W_{\b/2,\b/2} + 
\int_{-1}^{-\ge} \frac{dt}{2t} \\ \label{3}
-\int_{-1}^{-\ge} \frac{dt}{2t} - \int_{-\ge}^{u-\gl} \frac{dt}{2t} -
\int_{ u-\gl}^{-\gl} \frac{dt}{2t} .\end{align} 
The integrals in (\ref{3}) add up to $-\int_{ -1}^{-\gl} dt/2t$ and 
contribute $\gl^{-1/2}$ to the coefficients $d_{\pm}$. 
The sum in (\ref{2}) is a function of $\gl$
analytic near $\gl=u$ (since $\a'$ is $\sigma$-invariant) and is the same
for both cycles $\a_{\pm}$. The values of (\ref{1}) 
may depend on the cycle $\a_{\pm}$ but are independent of $\gl$.
We claim that in the limit $\ge \to 0$ the difference is an odd 
multiple of $\pi i$.
Indeed, transporting $\a_{-}$ along the composition
$\c_{-}\c_{+}^{-1}$ yields $\a_{+}$. On the other hand 
$2\W_{\a',\b}=\W_{\a_{+},\a_{+}}-\W_{\a_{-},\a_{-}}$. Thus the
difference of the two values of (\ref{1}) can be interpreted as
$\oint \W_{\a_{-},\a_{-}}$ along a loop $\c_{\ge}$ starting and terminating at
$\tau-(u+\ge)\1$ and transporting $\a_{-}$ to $\a_{+}$. Let us compose it
with a small loop $\d_{\ge}$ of radius $\ge$ around $\gl=u$. Since 
$\a_{+}$ transports along this loop back to $\a_{-}$, the composite integral
$\int_{\c_{\ge}\d_{\ge}}\W_{\a_{-},\a_{-}} \in 2\pi i \ZZ$ due to Proposition
$2$ and does not depend on $\ge$. Our claim follows therefore 
from Proposition $4$.

We conclude that $d_{\pm} = \pm d_0(\gl) \gl^{-1/2}$ where $d_0$ is 
a non-vanishing analytic function near $\gl=u$. Now we use the fact
that $\D_{A_1}$ is a tau-function of the KdV-hierarchy  
(\ref{KdV}) to conclude that the factor in (\ref{res2}) of the form 
\[ \sum_{\pm} \pm \frac{d\gl}{\sqrt{\gl}} \ \left[
 \ ^{A_1} \Gamma_{0}^{\pm\b/2}(\gl) \otimes \ ^{A_1}\Gamma_{0}^{\mp\b/2}(\gl)
\right] \ (\D_{A_1}\otimes \D_{A_1}) \]
is everywhere analytic in $\gl$. The same remains true after application
of the operator
$(\Psi \hat{R} e^{(U/z)\hat{\ }})^{\otimes 2}$. The vertex operator
$\Gamma_{\tau}^{\a'}\otimes \Gamma_{\tau}^{-\a'}$ is analytic near 
$\gl=u$ since $\a'$ is $\sigma$-invariant. Thus (\ref{res2}) has no pole
at $\gl=u$ and contributes $0$ to the residue sum.  

\medskip

Finally, consider the summands in (\ref{res}) with $\ga =\pm \b$.
Applying Theorems C and D, we transform the corresponding summands from
(\ref{res1}) to the form 
\[ \left(\Psi(\tau) \hat{R}_{\tau} e^{(U/z)\hat{\ }} 
\right)^{\otimes 2}
\left[  ... \D_{A_1}^{\otimes 2} \otimes 
\left( \sum_{\pm} e \gl d\gl   
\left( \ ^{A_1}\Gamma_{0}^{\pm\b} \otimes \ ^{A_1}\Gamma_{0}^{\mp\b}\right)  
\D_{A_1}^{\otimes 2} \right)^{(i)} \otimes \D_{A_1}^{\otimes 2}  ... \right] 
,\] 
where 
\[ e = \exp \left\{ -\int_{-1}^{-\ge}\frac{2dt}{t}-\int_{-\ge}^{u-\gl}
\frac{2dt}{t}-\int_{u-\gl}^{-\gl} \frac{2dt}{t} \right\} = 
\exp \left\{ -\int_{-1}^{-\gl} \frac{2dt}{t}\right\} = \gl^{-2} .\] 
The contribution of these terms to the residue sum (\ref{res1}) at $\gl=u^i$
can be calculated using the form (\ref{SL_2}) 
of the KdV-hierarchy for $\D_{A_1}$ and is equal to
\[ 16 \ (\Psi \hat{R} e^{(U/z)\hat{\ }})^{\otimes 2} \ 
l^{(i)} \ (\D_{A_1}^{\otimes N})^{\otimes 2}, \ \text{where} \ 
l^{(i)} =(... 1\otimes l \otimes 1 ...).\]   
In order to justify this conclusion, recall from the end of 
Section $5$ that conjugation
by $\exp (u^i/z)$ act as translation $\gl-u^i \mapsto \gl$. Also, since
$\D_{A_1}$ is tame, the vertex operator expression in (\ref{SL_2}) yields
a meromorphic $1$-form in $\gl$ with a singularity only at $\gl=0$. 
Thus the residue in (\ref{SL_2}) at $\gl=\infty$ is the same as at $\gl=0$.

Let us summarize our computation.

\medskip

{\bf Proposition 5.} {\em 
The residue sum (\ref{res}) is equal to
\begin{equation} \label{res3}  16\ 
\left( e^{c}\hat{S}^{-1} \Psi \hat{R} e^{(U/z)\hat{\ }} \right)^{\otimes 2} \ 
\left(\frac{N}{8}+\sum_{i=1}^N l^{(i)} \right) \ 
\left( \D_{A_1}^{\otimes N} \right)^{\otimes 2}.  
\end{equation}}

\medskip 

\subsection*{\bf 7. The Virasoro operator.} Functions of the form
$\Phi\otimes \Phi$ belong to a Fock space which is the quantization
of the symplectic space $\H\oplus \H$, the direct sum of two copies 
of $(\H,\Omega)$. Respectively the operator 
\begin{equation} \label{l_0} 
\sum_{k\geq 0} \sum_{a} (\frac{m_a}{h}+k) (q^a_k\otimes 1-1\otimes q^a_k)
(\frac{\p}{\p q^a_k}\otimes 1-1\otimes \frac{\p}{\p q^a_k}) \end{equation}
in (\ref{L-nod}) is
the quantization of a certain quadratic hamiltonian 
$\Omega (D \f, \f)/2$ on $\H\oplus \H$. Let us describe the 
infinitesimal symplectic transformation $D$ explicitly. 

Introduce the {\em Virasoro operator} $l_0:=z\p_z+1/2-\mu$. 
\footnote{The name comes from the property of the operators
$l_m:=l_0 z l_0 z ... z l_0$, ($z$ repeated $m$ times, $m=-1,0,1,2,...$)
to form a Lie algebra isomorphic to the algebra of formal vector fields 
$x^{m+1} \p/\p x$ on the line and participating in the formulation
of the Virasoro constraints (see \cite{GiQ, GiF}).}  
Since $\mu^*=-\mu$, the operator $l_0: \H \to \H$
is anti-symmetric with respect to $\Omega$, and the corresponding quadratic
hamiltonian reads
\[ \frac{1}{4\pi i}\oint (l_0 \f (-z) ,\f (z)) \ dz = \sum_{k\geq 0}
((k+\frac{1}{2}-\mu) f_k, (-1)^k f_{-1-k}) = 
-\sum_{k\geq 0}\sum_{a=1}^N (\frac{m_a}{h}+k)q_k^ap_k^a .\] 
Comparing this with (\ref{l_0}) we conclude that 
\begin{equation} \label{end} 
D = \left[ \begin{array}{rr} -l_0 & l_0 \\ l_0 & -l_0 \end{array} \right]
\ \in \  \operatorname{End} (\H \oplus \H) .\end{equation}
The expression (\ref{L-nod}) on the R.H.S. of the Hirota equation 
is proportional to
\[  \hat{D} \D^{\otimes 2} = \hat{M}^{\otimes 2} \ 
(\hat{M}^{\otimes 2})^{-1} \hat{D} (\hat{M}^{\otimes 2}) 
(\D_{A_1}^{\times N})^{\otimes 2}, \]
where $\hat{M} = e^{c(\tau)}\hat{S}_{\tau}^{-1}\Psi(\tau) \hat{R}_{\tau} 
e^{(U_{\tau}/z)\hat{\ }}$. Note that 
$\hat{M}^{\otimes 2}$ is the quantization of a block-diagonal
operator 
\[ B:=\left[ \begin{array}{cc} M & 0 \\ 0 & M \end{array} \right] ,\ 
\text{and} \ B^{-1} D B = \left[ \begin{array}{rr} -M^{-1} l_0 M & 
M^{-1} l_0 M \\ M^{-1} l_0 M & -M^{-1} l_0 M \end{array} \right] .\] 

\medskip

{\bf Proposition 6.} $M^{-1} l_0 M = \sum_{i=1}^N (\ ^{A_1}l_0)^{(i)}$.

\medskip

{\em Proof.} We have
\[ S (z\p_z+\frac{1}{2}-\mu) S^{-1} = z\p_z +\frac{1}{2}
-\mu + \frac{E\bullet}{z}, \] since 
$(z\p_z +L_E) S = \mu S-S \mu$ and  $\p_a S = z^{-1}\p_a\bullet S$. Next,
in the canonical coordinates $E=\sum u^i \p /\p u^i$, and therefore
\[ \Psi^{-1} (z\p_z +\frac{1}{2}-\mu + \frac{E\bullet}{z}) \Psi = 
z\p_z+\frac{1}{2}-V+\frac{U}{z}, \ \text{where}\ V:= \Psi^{-1}\mu \Psi = 
\Psi^{-1} L_E \Psi.\] Furthermore, the differential equations 
$\p_a (\Psi R e^{U/z}) = z^{-1}(\p_a \bullet) (\Psi R e^{U/z})$ translate
into $(d + \Psi^{-1}d\Psi) R = z^{-1} ( dU\ R-R\ dU)$. This implies 
$(L_E + V) R = z^{-1} (UR-RU)$, which together
with the homogeneity condition $(z\p_z +L_E) R = 0$ shows that
\[ R^{-1}(z\p_z+\frac{1}{2}-V+\frac{U}{z}) R = z\p_z+\frac{1}{2}+\frac{U}{z}.
\]
Finally 
\[ e^{-U/z} (z\p_z+\frac{1}{2}+\frac{U}{z}) e^{U/z} = z\p_z+\frac{1}{2}.
\hspace{2cm} \square \]

{\bf Proposition 7.} $\hat{M}^{-1} \hat{l}_0 \hat{M} = (M^{-1}l_0M)\hat{ }
+\tr \mu \mu^* /4 $.

\medskip

{\em Proof.} The quadratic hamiltonians for $z\p_z+1/2, \mu, V, \ln S, 
(E\bullet)/z, U/z$ contain no $p^2$-terms, and the quadratic hamiltonians for 
$z\p_z+1/2, \mu, V, \ln R$ contain no $q^2$-terms. Therefore, in the quantized
version of the previous computation, the only point where the 
cocycle $\C$ makes a non-trivial contribution is:
\[ \hat{R}^{-1}\ (\frac{U}{z})\hat{\ } \ \hat{R} = ( R^{-1} \frac{U}{z} R)
\hat{\ }\ +\ C.\]
Let $A=\ln R, B=U/z$. Then the quadratic hamiltonian of $BA-AB$ contains no 
$q^2$-terms (since $R|_{z=0}=1$). We have therefore
\begin{align} \notag 
\frac{d}{dt} e^{-t\hat{A}} \hat{B} e^{t\hat{A}} = e^{-t\hat{A}}
(\hat{B}\hat{A}-\hat{A}\hat{B})e^{t\hat{A}} = e^{-t\hat{A}}\ (BA-AB)\hat{\ }\ 
e^{t\hat{A}} + \C (B,A) =  \\ \notag
 \left[ e^{-tA} (BA-AB) e^{tA}\right]\hat{\ }\ +\ \C(B,A) =
\frac{d}{dt} \left[ e^{-tA} B e^{tA} \right]\hat{\ } + \C(B,A) .\end{align}
Integrating in $t$ from $0$ to $1$ we find $C = \C (B,A)$. Since 
$A=R_1z+o(z)$, we compute explicitly $C =\tr (BA)/2=\sum_i R_1^{ii} u^i/2$.  

This expression, which seems to be a function of $\tau$, has to be
a constant, and the value of this constant is well-known to be 
$\tr \mu \mu^* /4$ (see for instance the last chapter in \cite{H}).
For the sake of completeness we include the computation.
Namely, comparing the $z^0$- and $z^1$-terms in the equation 
$(L_E+V) R = z^{-1}(UR-RU)$ we find $V^{ij} =(u^i-u^j) R_1^{ij}$ and 
respectively 
\[ R_1^{ii}=-L_ER_1^{ii}=\sum_j (u^i-u^j)R_1^{ij}R_1^{ji} =
\sum_j \frac{V^{ij}V^{ji}}{u^j-u^i}.\]
Thus we have
\[ \frac{1}{2}\sum_i u^iR_1^{ii}=\sum_{ij} \frac{u^i V^{ij}V^{ji}}{2(u^j-u^i)}
= \sum_{ij} \frac{u^j V^{ij}V^{ji}}{2(u^i-u^j)} = 
-\frac{1}{4}\sum_{ij} V^{ij}V^{ji} = \frac{1}{4}\tr \mu \mu^* ,\]
since $V=\Psi^{-1}\mu \Psi$ and $V^t=\Psi^{-1}\mu^* \Psi$. $\square$    
 
\medskip

{\em Remark.} Slightly generalizing Propositions $6$ and $7$ one
obtains the following transformation formula (see Theorem $8.1$ in \cite{GiQ})
$\hat{M}^{-1} \hat{l}_m \hat{M} = \sum_i \ ^{A_1}\hat{l}_m^{(i)}$ 
for the Virasoro 
operators with $m\neq 0$. Since $(\hat{l}_m-\delta_{m,0}/16) \D_{A_1}=0$, 
this implies that $\D = \hat{M} \D_{A_1}^{\otimes N}$ satisfies the Virasoro
constraints $[\hat{l}_m -\delta_{m,0} \tr (\mu\mu^*/4+1/16)] \D =0$.
In fact this is Corollary $8.2$ in \cite{GiQ} specialized to Frobenius
structures of weighted-homogeneous singularities.

\medskip

Note that the conjugation $l_0 \mapsto M^{-1}l_0M$ of the {\em off}-diagonal
blocks in the matrix $D$ yields 
after quantization $\hat{M}^{-1}\hat{l}_0\hat{M} = (M^{-1}l_0M)\hat{\ }$ 
(since the cocycle $\C$ vanishes on pairs of quadratic 
hamiltonians corresponding to block-diagonal and block-off-diagonal 
operators.) Thus $\hat{B}^{-1} \hat{D} \hat{B} = \sum_i l^{(i)} 
- \tr \mu \mu^*/2$. Taking into account that  
\[ \frac{N (h+1)}{12 h}=\sum_a \frac{m_a (h-m_a)}{2h^2} = \frac{1}{2} \sum_a 
(\frac{1}{2}+\mu_a) (\frac{1}{2}-\mu_a) = 
\frac{1}{2}\tr (\frac{1}{4}+\mu\mu^*) \]
we conclude that the R.H.S. of the Hirota equation 
(\ref{Hirota},\ref{L-nod}) can be written as
\[ \hat{M}^{\otimes 2}\ (\frac{N}{8}+\sum_i l^{(i)})\  
(\D_{A_1}^{\otimes N})^{\otimes 2}.\]
Comparing this with Proposition $5$ we arrive at the following result.

\medskip

{\bf Proposition $8$.} {\em The function $\D$ satisfies the Hirota 
quadratic equation (\ref{Hirota},\ref{L-nod}) with $a_{\a}=b_{\a}/16$.}

\medskip

Since $\D\neq 0$, the Hirota equation is thus rendered consistent,
and the following corollary completes the proof of Theorem $1$.

\medskip

{\bf Corollary.} {\em The average value}
\[ \frac{1}{Nh}\sum_{\a\in A} b_{\a} = \frac{4(h+1)}{3h^2}. \]  

\medskip

Note that in the proof of Theorem $1$ we use neither 
the Virasoro constraints for 
$\D_{A_1}$ nor the fact that the $N$ factors in $\D_{A_1}^{\otimes N}$ 
are the same. The only relevant conditions for $\D_{A_1}$ were both forms
of the $KdV$-hierarchy and the tame property of $\exp (u/z)\hat{\ } \D_{A_1}$.
Thus we have actually proved the following generalization of Theorem $1$.
 
\medskip

{\bf Theorem 2.} {\em Suppose that tame asymptotical functions
$\Phi_1, ..., \Phi_N$ are tau-functions of the KdV-hierarchy and
remain tame under the string flow 
$\Phi_i \mapsto \exp (u/z)\hat{\ } \Phi_i$ for all $u$. Then 
\[ \Phi :=e^{c(\tau)}\ \hat{S}_{\tau}^{-1}\ \Psi(\tau)\ \hat{R}_{\tau}\ 
\ e^{(U_{\tau}/z)\hat{\ }}\ (\Phi_1 \otimes ... \otimes \Phi_N) \]
satisfies the corresponding Hirota quadratic equation 
(\ref{Hirota} -- \ref{c}).}

\medskip  

{\em Remark.} Although the condition for $\Phi_i$ to remain tame under 
the string flow is quite restrictive, $\D_{A_1}$ is not the
only tau-function satisfying it. A large class of examples consists of the
shifts $\D_{A_1} (\q + {\mathbf a})$ where ${\mathbf a} (z) 
=a_0+a_1z+a_2z^2+...$
is a series with coefficients $a_k$ which are arbitrary series in $\h$ 
such that
$a_0$ and $a_1$ are smaller than $1$ in the $\h$-adic norm and 
$a_k \to 0 $ in this norm as $k \to \infty$.      

\medskip
  
\subsection*{\bf 8. The Kac -- Wakimoto hierarchies.} Let us compare
the ADE- hierarchies (\ref{Hirota}--\ref{c}) with the 
{\em principal hierarchies}
of the types $A_N^{(1)}, D_N^{(1)}, E_N^{(1)}$ described in Theorem $1.1$
in \cite{KW}. The corresponding Hirota equation $(1.14)$ in \cite{KW} 
has the form 
\begin{align} \notag \Res \frac{d\zeta}{\zeta} \sum_{i=1}^N g_i 
e^{\sum_{m\in E_{+}} 2\b_{i,\overline{m}} \h^{-1/2}y_{m} \zeta^m} 
e^{-\sum_{m\in E_{+}} \b_{i,-\overline{m}} \h^{1/2}\p_{y_m} \zeta^{-m}/m} 
\Phi (\x+\y) \Phi (\x-\y) \\ \label{KW} = \ 
\left( 2h \sum_{m\in E_{+}} m\ y_m  \ + 
\lan \rho,\rho\ran \right) \Phi (\x+\y) \Phi(\x-\y).  
\end{align}   
Here $\rho$ is the sum of the fundamental weights of the root system $A$,
and the value $\lan \rho ,\rho \ran = N h (h+1) /12$ can be found 
for instance from the tables in \cite{Bo}. The index set 
$E_{+} = \{ m_a+kh | a=1,...,N, \ k=0,1,2,...\}$, and 
$\overline{m}$ denote the remainder modulo $h$. The vertex operators in the
sum correspond to a set of roots $\a_i, i=1,...,N$, chosen one from 
each orbit of some Coxeter element $M$ on the root system $A$. The 
coefficients $\b_{i,\overline{m}}$ are coordinates of $\a_i$ with respect
to a basis of eigenvectors $H_{\overline{m}}$ of the Coxeter transformation
$M$ with the eigenvalues $\exp (2\pi \sqrt{-1} m/h)$. The coefficients $g_i$
are defined via representation theory of affine Lie algebras. 
The numerical values of $g_i$ are computed in \cite{KW} in the cases 
$A_N$, $D_4$ and $E_6$.
 
In order to identify the vertex operators in (\ref{KW}) with those in
(\ref{Hirota},\ref{L-nod}) let us start with taking 
$\zeta = (h\gl)^{1/h}$. Then the components of the period vector 
$I^{(-1)}_{\a_i}$ with respect to a suitable basis $[\psi_a]\in H$ will
have the form
\begin{equation} \label{psi} 
(I^{(-1)}_{\a_i}(\gl), [\psi_a])=\b_{i,\overline{m_a}} \ m_a^{-1}\ 
(h\gl)^{m_a/h} \end{equation}
since the weighted - homogeneous forms $\psi_a \go /df$ represent 
a basis of eigenvectors for the classical monodromy operator in 
$H^2(f^{-1}(1),\CC)$. Then it is straightforward to check that the
relation 
\begin{equation} \label{x} q_k^a = \prod_{r=0}^k (m_a+rh) \ t_{m_a+kh} 
\end {equation} 
(together with the standard change $\x+\y=\t',\ \x-\y =\t''$ as in Section
$4$) identifies the vertex operators in (\ref{KW}) with 
$\Gamma^{\a_i}\otimes \Gamma^{-\a_i}$.
Note that replacing $\a_i$ with any of the $h$ roots from the same $M$-orbit
does not change the corresponding residue in (\ref{KW}) since 
the new vertex operator would differ from the old one only by the choice
of the branch of $\zeta = (h\gl)^{1/h}$. Thus we arrive at the following
conclusion.

\medskip     

{\bf Proposition 9.} {\em The choice of the basis $\{ [\psi_a]\in H \}$ such
that (\ref{psi}) holds true and the change of variables (\ref{x})
identify the Hirota equation (\ref{Hirota}--\ref{c}) with the corresponding
hierarchy of the form (\ref{KW}) provided that $g_i = h^3 a_{\a_i}=
h^3 b_{\a_i}/16$.}

\medskip

Let us now compute the coefficients $b_{\a}$. First, rewrite the definition
(\ref{b}) as 
\[ b_{\a} = \lim_{\ge \to 0} 
e^{-\int_{-\ge \1}^{\tau -(u+\ge)\1} \W_{\a,\a}} =
\lim_{\ge \to 0} \prod_{\c\in A} 
\frac{\lan y(\ge),\c \ran^{\frac{\lan \a,\c\ran^2}{2}}}
{\lan \ge^{1/h}\kappa, \c \ran^{\frac{\lan \a,\c\ran^2}{2}}} = 
\frac{\lan v,\a\ran^4}{\lan \kappa,\a \ran^4} \prod_{\lan \c,\a \ran=1}
\frac{\lan x,\c\ran}{\lan \kappa,\c\ran} ,\] 
where $\ge^{1/h}\kappa$, $y(\ge)$ and $x$ are inverse images under the 
Chevalley map of $-\ge \1$, $\tau-(u+\ge)\1$ and $\tau-u\1$
respectively, $x$ is a generic point on the mirror $\lan \a,x\ran = 0$
and $v$ is determined from the expansion $y(\ge) = x+\ge^{1/2}v+o(\ge^{1/2})$.
We will use this formula in the case of $A$ and $D$ series.

\medskip

{\em Case $A_N$.} The root system consists of the vectors 
$\c_{ij}:=e_i-e_j$ in the space $\CC^{N+1}$ with the standard orthonormal 
basis $e_0,...,e_N$ and coordinates $z_0,...,z_N$. Take 
\[ F(z,\tau) =\frac{z^{N+1}}{N+1}+t_1z^{N-1}+...+t_N = 
\frac{1}{N+1} \prod_{i=0}^N (z-z_i).\]
Let $\a=e_a-e_b$ and let $t=\tau-u\1$ be a generic point on the discriminant.
Then the components $y_i(\ge) = x_i+\ge^{1/2}v_i+\ge w_i + o(\ge)$
(where $x_a=x_b$) satisfy $F(y_i,t-\ge \1) = 0$ and therefore
\[ \ge = F(x_i,t)+F'(x_i,t) \ge^{1/2}y_i+F'(x_i,t)\ge w_i +
F''(x_i,t)\ge \frac{v_i^2}{2}+o(\ge) .\]
We have $F(x_i,t)=0$ for all $i$ and $F'(x_i,t)=0$ for $i=a,b$. This implies
that $v_i=\pm \sqrt{2/F''(x_a,t)}$ for $i=a,b$ and hence 
$\lan \a, v \ran  = \pm 2\sqrt{2/F''(x_a,t)}$. Thus $\lan \a,v\ran^4 =
64/ F''(x_a,t)^2$. On the other hand,
\[ \prod_{\lan \c,\a\ran =1} \lan x,\c\ran = (-1)^{N-1}
\prod_{i\neq a,b} (x_i-x_a)^2 = (-1)^{N-1}\left( 
\frac{(N+1)F''(x_a,t)}{2} \right)^2 .\]
The eigenvector $\kappa = (N+1)^{1/(N+1)}(1,\eta,\eta^2,...,\eta^{N-1})$ 
of the Coxeter transformation $(z_0,...,z_N)
\mapsto (z_1,...,z_N,z_0)$ with the eigenvalue $\eta = \exp 2\pi i/(N+1)$ is
a preimage of $t =-\1$ under the Chevalley map.  
We find \footnote{We use here the facts that the product 
$\prod_{k\neq a} (\zeta-e^{2\pi i k/n})$ over all $n$-th roots of unity 
except $\zeta = e^{2\pi i a/n}$ is equal to the derivative of $z^n-1$ at 
$z=\zeta$, {\em i.e.} to $n/\zeta$.}
\begin{align} \notag
 \lan \kappa,\a \ran^4\prod_{\lan \a,\c\ran=1}\lan \kappa,\c\ran =
(N+1)^2 (\eta^a-\eta^b)^2 \prod_{j\neq a}(\eta^a-\eta^j)\prod_{i\neq b}
(\eta^i-\eta^b) = \\ \notag (-1)^N(N+1)^4 (\eta^a-\eta^b)^2 \eta^{N(a+b)}
 = (-1)^{N-1}(N+1)^4 (2-\eta^{a-b}-\eta^{b-a}).\end{align}
Collecting the results we find
\[ b_{\a} = \frac{16}{(N+1)^2}\frac{1}{(2-\eta^{a-b}-\eta^{b-a})} .\]
This agrees with Theorem $1.2$ in \cite{KW} where 
$g_i = (N+1)/(2-\eta^i-\eta^{-i})$ corresponds to $\a_i = e_0-e_i$.
In particular 
\[ \sum g_k = \frac{N+1}{4} \sum_{k=1}^N \sin^{-2} (\frac{\pi k}{N+1}) =
\frac{N(N+1)(N+2)}{12} .\]
The middle expression is a special case of {\em Dedekind sums}, and the
second equality, which follows from our results, is well known in number
theory (see {\em i.g.} \cite{Be}).

\medskip

{\em Case $D_N$.} The root system consists of the vectors
$ \pm e_i \pm e_j\ , \ i\neq j$, where $e_1, ... e_N$ is the standard 
orthonormal basis and $(z_1,...,z_N)$ are the corresponding coordinates
in $\CC^N$. The parameters $(t_1,...,t_N)$ in the following family of 
polynomials
\[ F(z,t) = z^{2N}+t_2z^{2N-2}+t_3z^{2N-4}+...+t_Nz^2 + t_1^2 = \prod_{i=1}^N
(z^2-z^2_i) \]
are identified with coordinates on the Chevalley quotient $\CC^N/W$.
Note that the invariant $t_N$ of degree $h=2N-2$ is the coefficient at $z^2$.
Let us assume that $x$ is a generic point on the mirror $z_a\pm z_b$ 
orthogonal to the root $\a = e_a \mp e_b$, and that $t$ 
is the corresponding point on the discriminant, so that $x_a=\pm x_b$ and
$F(\pm x_a,t)=F'(\pm x_a,t)=0$. Taking 
$y_i(\ge)=x_i+\ge^{1/2}v_i+\ge w_i + o(\ge)$ and 
expanding $F(y(\ge),t_1,...,t_N-1,t_N-\ge)=0$ in $\ge$ we find
\[ \ge x_i^2 = F(x_i,t) + F'(x_i,t) (\ge^{1/2}v_i+\ge w_i) + 
F''(x_i,t) \ge \frac{v_i^2}{2} + o (\ge) .\]
Thus $v_a = \sqrt{2 x_a^2/ F''(x_a,t)}$, $v_b=\mp \sqrt{2 x_a^2/F''(x_a,t)}$
and $\lan \a, v\ran^4 = 64 x_a^4 / F''(x_a,t)^2$. Furthermore,
\[ F''(z,t)|_{z=x_a}=[
2 \sum_a \prod_{i\neq a}(z^2-x_i^2) + 4z^2\sum_{a\neq b}
\prod_{i\neq a,b} (z^2-x_i^2) ]_{z=x_a}=8 x_a^2 
\prod_{i\neq a} (x_a^2-x_i^2).\]
Using this we find     
\[ \prod_{\lan \a,\c\ran =1} \lan x,\a \ran = 
\prod_{j\neq a,b} (x_a^2-x_j^2) \prod_{i\neq a,b} (\mp 1)(x_i^2-x_b^2) =  
(\pm 1)^{N-2} \frac{F''(x_a,t)^2}{64 x_a^4}.\]
Next, the eigenvector $\kappa = (1,\eta,...,\eta^{N-2},0)$ of 
the Coxeter transformation 
\[ (z_1,...,z_N) \mapsto (z_2,...,z_{N-1},-z_1,z_N) \]
with the eigenvalue $\eta = \exp \pi i /(N-1)$ is mapped
to $(t_1,...,t_N)=(0,...,0, 1)$ under the Chevalley map.
Assuming first that $\a=e_a\mp e_b$ with $a,b<N$ we find
\begin{align} \notag  
\lan \kappa, \a\ran^4\prod_{\lan \a,\c\ran=1}\lan \kappa, \c\ran  =  
(\eta^a\mp \eta^b)^4 (\eta^a)^2(\mp \eta^b)^2 
\prod_{i\neq a,b,N} (\eta^{2a}-\eta^{2i}) (\pm 1)
(\eta^{2b}-\eta^{2i}) = \\ \notag -(\pm 1)^{N-2} (N-1)^2
\frac{(\eta^a\mp \eta^b)^2}{(\eta^a\pm \eta^b)^2} = (\pm 1)^{N-2} (N-1)^2
\frac{(2\mp \eta^{a-b}\mp \eta^{b-a})}{(2\pm \eta^{a-b}\pm\eta^{b-a})}.
\end{align}
Combining with the previous formulas we compute
\[ b_{\a} = \frac{1}{(N-1)^2} 
\frac{(2\pm \eta^{a-b}\pm \eta^{b-a})}{(2\mp \eta^{a-b}\mp\eta^{b-a})} 
\ \ \text{for}\ \ \a =e_a\mp e_b.\]
Now let $\a =e_a\mp e_N$. Then  
\[ \lan \kappa, \a\ran^4\prod_{\lan \a,\c\ran=1}\lan \kappa, \c\ran = 
\eta^{4a} \prod_{i\neq a,N} (\eta^{2a}-\eta^{2i}) (-\eta^{2i}) =  
%\eta^{2a} \eta^{-2a}(N-1) (-1)^{N-2} \eta^{(N-2)(N-1)-2a} =
(-1)^{N-2}(N-1)\]
and therefore $b_{\a} = 1/(N-1)$.

Taking the representatives 
\[ \a_1=e_{N-1}-e_1,\ ...,\ \a_{N_2}=e_{N-1}-e_{N-2},\ 
\a_{N-1}=e_{N-1}-e_N,\ \a_N= e_{N-1}+e_N \]
in the orbits of the Coxeter transformation on $A$, we find
\[ g_i = \frac{(N-1)}{2} \frac{(2-\eta^i-\eta^{-i})}{(2+\eta^i+\eta^{-i})}
\ \text{for} \ i=1,...,N-2,\ \text{and}\ g_i=\frac{(N-1)^2}{2} \ \text{for}\
i=N-1,N.\]
The identity $\sum g_k = (N-1)N(2N-1)/6$, which follows from our general 
theory, agrees with the value of the Dedekind sum
\footnote{It is essentially the same one as in the $A$-case since 
$\sin^{-2}x - 1 =\cot^2x =\tan^2(\pi/2-x)$.}
 \[ \sum_{k=1}^{N-2}\tan^2 (\frac{\pi k}{2N-2})  = \frac{(N-2)(2N-3)}{3} .\]

In the case $N=4$ the values $g_i=1/2, 9/2, 9/2, 9/2$ agree with the values 
of $g_i$ found in \cite{KW}, Proposition $1.3 (a)$. 

\medskip

{\em Cases $E_N$.} We find $g_i$ using the packages {\em LiE} and {\em MAPLE} 
to compute the ratios via (\ref{c}) and then apply the normalizing relation
(\ref{c0}). In each case $E_N$, let $\a_1,...,\a_N$ be the simple roots
and $M$ be the Coxeter transformation described by the following diagrams:
\[ \begin{array}{l} \begin{array}{ccccccccl}
  \a_1 & &  \a_3 & & \a_4  & & \a_5  & & \a_6   \\
\bullet&-&\bullet&-&\bullet&-&\bullet&-&\bullet \\
       & &       & &   |   & &       & &        \\
       & &       & &\bullet& &       & & M = \gs_1\gs_4\gs_6\gs_2\gs_3\gs_5,\\
       & &       & & \a_2  & &       & & \end{array} \\ 
 \begin{array}{ccccccccccl}
  \a_1 & &  \a_3 & & \a_4  & & \a_5  & & \a_6  & &  \a_7  \\
\bullet&-&\bullet&-&\bullet&-&\bullet&-&\bullet&-&\bullet \\
       & &       & &   |   & &       & &       & &        \\
       & &       & &\bullet& &       & &       & & 
           M = \gs_1\gs_4\gs_6\gs_2\gs_3\gs_5\gs_7,       \\
       & &       & & \a_2  & &       & &       & &        \end{array} \\ 
 \begin{array}{ccccccccccccl}
  \a_1 & &  \a_3 & & \a_4  & & \a_5  & & \a_6  & &  \a_7 & &  \a_8  \\
\bullet&-&\bullet&-&\bullet&-&\bullet&-&\bullet&-&\bullet&-&\bullet \\
       & &       & &   |   & &       & &       & &       & &        \\
       & &       & &\bullet& &       & &       & &       & & 
               M = \gs_1\gs_4\gs_6\gs_8 \gs_2\gs_3\gs_5\gs_7.       \\
       & &       & & \a_2  & &       & &       & &       & &        \end{array}
 \end{array} \]
One can check ({\em i.g.} using {\em LiE}) that all simple roots 
$\a_1,...,\a_N$ belong to different $M$-orbits. 
The following tables represent the values of the corresponding 
coefficients $g_i$ while the values of $b_{\a_i}$ can be obtained
from them as in Proposition $9$.   
 
\medskip

{\em Case $E_6$.} We have $b_{\a_i}=g_i/108$, where
\[ g_1=g_6=16+8\sqrt{3},\ 
g_3=g_5=16-8\sqrt{3},\  g_2=7+4\sqrt{3}, \ g_4=7-4\sqrt{3}.\]
This agrees with the values of $g_i$ found
in \cite{KW}, Proposition $1.3(b)$.

\medskip

{\em Case $E_7$.} We have $b_{\a_i}=2 g_i/729$. Put $u=\cos (\pi/9)$. Then
\[ g_1=\frac{27}{2}+36u+24u^2, \  
g_2=\frac{225}{2}+36u-144u^2, \ g_3=\frac{3}{2}, \  
g_4=\frac{147}{2}+12u-96u^2, \]
\[  g_5=\frac{9}{2}-72u+72u^2, \ \  
g_6=-\frac{21}{2}-48u+72u^2, \ \ g_7=\frac{9}{2}+36u+72u^2 \ .\]

\medskip

{\em Case $E_8$.} We have $b_{\a_i}=2g_i/3375$. Let $u=\cos (\pi/15)$. Then
\begin{align} \notag g_1=\frac{33}{2}+80u+72u^2-16u^3,\ \  &  
   g_2=\frac{273}{2}+132u-136u^2-128u^3, \\ \notag
   g_3=-\frac{123}{2}+568u+376u^2-912u^3,\ \  &  
   g_4=\frac{109}{2}-368u-72u^2+400u^3,\\ \notag
   g_5=\frac{745}{2}+584u-376u^2-624u^3,\ \  &  
   g_6=\frac{257}{2}-1220u-232u^2+1376u^3,\\ \notag
   g_7=-\frac{35}{2}+156u+136u^2-256u^3,\ \  & 
   g_8=-\frac{19}{2}+68u+232u^2+160u^3\ .\end{align}

\medskip

\subsection*{\bf 9. Open questions.} 

(a) The formula (\ref{D}) defines
the total descendent potential $\D$ as an asymptotical function of 
$\q=\q_0+\q_1z+...$ with semisimple $\q_0$. As it is shown in \cite{GiI},
Theorem $5$, the function $\D_{A_N}$ extends to arbitrary values 
of $\q_0$ without singularities. We expect the same for $\D_{D_N}$ and
$\D_{E_N}$ but leave this issue open.

(b) B. Dubrovin \cite{D} associates to a Frobenius manifold a 
{\em dispersionless} integrable hierarchy. In particular, the hierarchy 
(\ref{Hirota},\ref{L-nod}) for asymptotical
functions $\Phi = \exp (\F^{(0)}/\h+\F^{(1)}+...)$ admits the
{\em dispersionless limit} as $\h\to 0$ which is an infinite
system of equations for $\F^{(0)}$. It is not hard to show that
{\em $\F $ satisfies the dispersionless hierarchy if and only if 
the Gaussian distributions $\Phi :=\exp \{d^2_{\x}\F (\q)/2\h\}$
(where $d^2_{\x}\F$ is the quadratic differential of $\F$ at $\x$) 
satisfy the original hierarchy (\ref{Hirota},\ref{L-nod}) for all $\x$.}
An elegant explicit characterization in terms of the semi-infinite
Grassmannian of those Gaussian distributions which satisfy
the hierarchy of the type $A_N$ is given in the appendix to \cite{GiI}.
It would be interesting to generalize the characterization to the cases
$D_N$, $E_N$.

(c) Theorem $1$ implies that the genus $0$ descendent potential
$\F^{(0)} =\lim_{\h\to 0} \h \ln \D$ satisfies the corresponding
dispersionless hierarchy. The quadratic forms $d^2_{\x}\F^{(0)}(\q)$ 
depend only on $N$ parameters $\tau =\tau(\x)$ (due to the property
$(\ast)$ of the cone $\L = \operatorname{graph} d\F^{(0)}$) and have
the following explicit description (see Appendix in \cite{GiI}):
\[ \int_{\0}^{\tau} \sum_a ( [S_{t}(z)\q(z)]_0\bullet [S_{t}(z)\q(z)]_0, 
\p_{t^a}) \ dt^a , \]
where $[S (z) \q(z)]_0=S_0\q_0+S_1\q_1+...$ denotes the $z^0$-mode. The
corresponding Gaussian distributions satisfy therefore the hierarchy 
(\ref{Hirota} -- \ref{c}). Taking $\tau = u\1$ so that 
$[S_{\tau}\q ]_0 =\sum \q_k u^k/k!$ we conclude that in particular 
the hierarchy has the $1$-parametric family of Gaussian solutions
\[ \Phi = \exp \left\{ \frac{1}{2\h} 
\int_0^u (\sum_{k\geq 0} \q_k\frac{v^k}{k!}, \sum_{l\geq 0} \q_l\frac{v^l}{l!})
\ dv \ \right\} .\]
This imposes non-trivial constraints on the coefficients $a_{\a}$ in
the Hirota equation (\ref{Hirota}). It would be interesting to find out if
these constraints are sufficient in order to determine the coefficients
unambiguously.

(d) Our computations in Section $8$ confirm the Conjecture from 
Section $1$ in the cases $A_N, D_4, E_6$ and leave it open in the cases
$D_N$ with $N>4$, $E_7$ and $E_8$ --- mostly because the values of the
coefficients $g_i$ in the Kac -- Wakimoto theory remain unknown. 
A more conceptual approach to the identification of the Hirota equations
should rely on the definition of the coefficients $g_i$ given in 
\cite{KW} in terms of representation theory. Namely, the vertex operators 
$C_i^{\pm} \Gamma^{\pm \a_i}$ participate in the so called {\em 
principal} construction of the basis representation of the affine Lie
algebra $\hat{A}_N, \hat{D}_N$ or $\hat{E}_N$, and $g_i=C_i^{+}C_i^{-}$.
Here $C_i^{\pm}$ are certain structure constants whose values remain 
generally speaking unknown. Our successful description of the 
products $C_i^{+}C_i^{-}$ via the phase forms suggests that one should
look for the intrinsic role of the phase forms in representation theory 
and for a description of the individual coefficients $C_i^{\pm}$ in terms
of the phase forms or their generalizations.

(e) In representation theory, the hierarchies of the ADE-type form only 
a part of a larger list of examples including twisted versions
of the affine Lie algebras and non-simply laced Dynkin diagrams.
It would be interesting to find the corresponding constructions
in singularity theory and, in particular, to associate the Hirota
equations to the boundary singularities $B_N, C_N, F_4$.

(f) B. Dubrovin and E. Zhang \cite{DZ} associate an integrable hierarchy
to any semisimple Frobenius manifold. In a sense their construction
is parallel to the definition (\ref{D}) of the total descendent potential
$\D$ (see \cite{GiQ}) and in particular yields objects defined
in the complement to the caustic. In this regard the vertex operator
description of the hierarchies seems more attractive as it is free of this
defect. Of course, the ADE-hierarchies (\ref{Hirota} -- \ref{c}) are
expected to be equivalent to the hierarchies of Dubrovin -- Zhang.
It would be interesting to confirm this expectation.

(g) Conjecturally, the total descendent potential $\D$ extends 
analytically across the caustic values of $\q_0$ in the case of
K. Saito's (semisimple) Frobenius structure corresponding to any
isolated singularity. (By the way, this is known to be false for,
say, boundary singularities or for finite reflection groups other than
$A_N, D_N, E_N$.) Respectively, one should expect the same for
the hierarchies of Dubrovin -- Zhang. It would be very interesting
to give a vertex operator description of the hierarchies
together with Theorem $1$ for arbitrary (or at least weighted - homogeneous) 
isolated singularities of functions. The most obvious difficulty 
is that the vertex operator sum (\ref{Hirota}) over the set of all
vanishing cycles (or even orbits of the classical monodromy operator on 
this set) becomes infinite beyond the $ADE$ list. Nevertheless we believe that
the obstructions can be removed by an appropriate generalization 
of the concepts involved. The first examples to study here would be 
the unimodal singularities $P_8, X_9, J_{10}$ (see \cite{Ar}). Their miniversal
deformations are closely related to the complex crystallographic reflection
groups $\tilde{E}_6, \tilde{E}_7, \tilde{E}_8$ (see \cite{Lo}). Moreover,
the question can be extrapolated to the complex crystallographic groups
$\tilde{A}_N,\tilde{D}_N$, and the $3$-dimensional Frobenius manifold
to be called $\tilde{A}_1$ represents the first challenge.

\medskip

\subsection*{Acknowledgments.} The authors are thankful to
V. Kac and M. Wakimoto for very helpful consultations on the material of 
their paper \cite{KW}, to E. Frenkel and P. Pribik for their 
interest and stimulating discussions and to B. Sturmfels for his 
recommendation of the {\em LiE} package.


\newpage

\begin{thebibliography}{1000}

\bibitem{Ar} V. I. Arnold. {\em Normal forms of functions near degenerate
critical points, Weyl groups $A_k, D_k, E_k$ and Lagrangian singularities.}
Functional Anal. Appl. {\bf 6} (1972), no. 2, 3 -- 25.

\bibitem{AVG} V. I. Arnold, S. M.  Gousein-Zade, A. N. Varchenko. 
{\em Singularities of differentiable maps. Vol. II. 
Monodromy and asymptotics of integrals.} 
Monographs in Mathematics, 83. Birkh\"auser Boston, Inc., Boston, MA, 1988. 
viii+492 pp

\bibitem{Ba} S. Barannikov. {\em Quantum periods - I. Semi-infinite variations
of Hodge structures.} Internat. Math. Res. Notices 2001, no. 23, 1243--1264. 
alg-geom/0006193.

\bibitem{Be} M. Beck. {\em The reciprocity law for Dedekind sums via the 
constant Ehrhart coefficient.} arXiv: math.NT/0305404. 

\bibitem{Bo} N. Bourbaki. {\em Groupes et alg\`ebras de Lie.} Ch. IV -- VI,
Hermann, $1968$. 

\bibitem{D} B. Dubrovin. {\em Geometry of 2D topological filed theories.} In:
Integrable Systems and Quantum Groups. 
Springer Lecture Notes in Math. 1620 (1996), 120--348. 

\bibitem{D2} B. Dubrovin. {\em Painlev\'e transcendents in two-dimensional 
topological field theory.} The Painlev\'e property, 287--412,
CRM Ser. Math. Phys., Springer, New York, 1999. arXiv: math.AG/9803107

\bibitem{DZ} B. Dubrovin, Y. Zhang. {\em Normal forms of hierarchies of 
integrable PDEs, Frobenius manifolds and Gromov -- Witten invariants.} 
arXiv: math.DG/0108160.

\bibitem{GiQ} A. Givental. {\em Gromov -- Witten invariants
and quantization of quadratic hamiltonians.} Moscow Mathematical Journal, 
v.1(2001), no. 4, 551--568.

\bibitem{GiI} A. Givental. 
{\em $A_{n-1}$-singularities and $n$KdV hierarchies.}
To appear in Moscow Mathematical Journal. arXiv: math.AG/0209205

\bibitem{GiA} A. Givental. {\em Asymptotics of the intersection form of a 
quasi-homogeneous function singularity.} Funct. Anal. Appl. 16 (4), 1982, 
294-297.

\bibitem{GiF} A. Givental. {\em Symplectic geometry of Frobenius structures.}
To appear in Proceedings of the workshop on Frobenius 
structures held at MPIM Bonn in July 2002. arXiv: math.AG/0305409

\bibitem{H} C. Hertling. {\em Frobenius manifolds and moduli spaces for 
singularities.}
Cambridge Tracts in Mathematics. Cambridge University Press, 2002, 280 pp. 

\bibitem{K} V. Kac. {\em Infinite dimensional Lie algebras.} 3rd edition. 
Cambridge University Press, 1990, xxii+400 pp.

\bibitem{KW} V. Kac, M. Wakimoto. {\em Exceptional hierarchies of soliton
equations.} Proceedings of the 1987 conference on theta functions in Maine,
Proceedings of Symposia in Pure Math. {\bf 49}, 1989, 138 -- 177. 

\bibitem{Ko}  M. Kontsevich. {\em Intersection theory on the moduli space of
curves and the matrix Airy function.} Commun. Math. Phys. {\bf 147} 
(1992), 1 -- 23.

\bibitem{Lo} E. Looijenga. {\em On the semiuniversal deformation of 
a simple-elliptic hypersurface singularity.} Topology {\bf 17} (1977), 
23 -- 40.

\bibitem{S} K. Saito. {\em On a linear structure of the quotient variety 
by a finite reflection group.} Publ. Res. Inst. Math. Sci. 29 (1993),
no. 4, 535--579. 

\bibitem{W} E. Witten. {\em Two-dimensional gravity and intersection theory
on moduli space.} Surveys in Diff. Geom.  1 (1991), 243--310.

\end{thebibliography}
\enddocument